\documentclass[a4paper]{article}
\usepackage{amssymb}
\usepackage{amsfonts}
\usepackage{amsmath}
\usepackage{amsthm}
\usepackage{amsxtra}
\usepackage{amscd}

\numberwithin{equation}{section}

{\theoremstyle{definition}\newtheorem{definition}{Definition}[section]
\newtheorem{notation}[definition]{Notation}
\newtheorem{terminology}[definition]{Terminology}
\newtheorem{remark}[definition]{Remark}
}
\newtheorem{proposition}[definition]{Proposition}
\newtheorem{lemma}[definition]{Lemma}
\newtheorem{theorem}[definition]{Theorem}

\newcommand{\cG}{\mathbb{G}}
\newcommand{\de}{\Delta}
\newcommand{\ot}{\otimes}

\newcommand{\cGh}{\widehat{\mathbb{G}}}

\newcommand{\recht}{\rightarrow}
\newcommand{\B}{\operatorname{B}}
\newcommand{\oH}{\overline{H}}

\newcommand{\sde}{\delta}
\newcommand{\vphi}{\varphi}
\newcommand{\Mor}{\operatorname{Mor}}
\newcommand{\Hom}{\operatorname{Hom}}

\newcommand{\ox}{\overline{x}}

\newcommand{\eps}{\epsilon}
\newcommand{\veps}{\varepsilon}
\newcommand{\dimq}{\operatorname{dim}_q}

\newcommand{\om}{\omega}

\newcommand{\Tr}{\operatorname{Tr}}
\newcommand{\multq}{\operatorname{mult}_q}
\newcommand{\mult}{\operatorname{mult}}

\newcommand{\R}{\mathbb{R}}
\newcommand{\C}{\mathbb{C}}
\newcommand{\la}{\langle}
\newcommand{\ra}{\rangle}
\newcommand{\cB}{\mathcal{B}}

\newcommand{\cC}{\mathcal{C}}
\newcommand{\cst}{C$^*$}

\newcommand{\deh}{\hat{\Delta}}
\newcommand{\N}{\mathbb{N}}

\newcommand{\al}{\alpha}

\newcommand{\GL}{\operatorname{GL}}

\newcommand{\meq}{\underset{\text{\rm mon}}{\sim}}
\newcommand{\SU}{\operatorname{SU}}
\newcommand{\SO}{\operatorname{SO}}

\newcommand{\M}{\operatorname{M}}
\newcommand{\U}{\operatorname{U}}

\newcommand{\cD}{\mathcal{D}}
\newcommand{\cL}{\mathcal{L}}

\newcommand{\tb}{\tilde{b}}
\newcommand{\id}{\mathalpha{\text{\rm id}}}
\newcommand{\fancyot}{\mathbin{\text{\footnotesize\textcircled{\tiny \sf T}}}}
\newcommand{\Irred}{\operatorname{Irred}}

\newcommand{\hon}{H^{\infty}(\widehat{\mathbb{G}},\mu)}
\newcommand{\hone}{H^{\infty}(\widehat{\mathbb{G}}_1,\mu)}
\newcommand{\hont}{H^{\infty}(\widehat{\mathbb{G}}_2,\mu)}
\newcommand{\plon}{\ell^{\infty}(\widehat{\mathbb{G}})}
\newcommand{\plone}{\ell^{\infty}(\widehat{\mathbb{G}}_1)}
\newcommand{\plont}{\ell^{\infty}(\widehat{\mathbb{G}}_2)}

\newcommand{\mue}{\Vb_1}
\newcommand{\mut}{\Vb_2}
\newcommand{\epsh}{\widehat{\varepsilon}}

\newcommand{\Vb}{\mathbb{V}}
\newcommand{\Xb}{\mathbb{X}}
\newcommand{\Yb}{\mathbb{Y}}

\newcommand{\wlim}[1]{w^*{\text{\rm -}}\lim_{#1}}

\newcommand{\cH}{\mathbb{H}}
\accentedsymbol{\Ghh}{\,\widehat{\!\widehat G}}

\newcommand{\diag}{\operatorname{diag}}
\newcommand{\End}{\operatorname{End}}

\newcommand{\Ad}{\operatorname{Ad}}

\newcommand{\otalg}{\ot_{\operatorname{alg}}}
\newcommand{\boxalg}{\boxtimes_{\operatorname{alg}}}
\newcommand{\otvn}{\overline{\ot}}
\newcommand{\omu}{\overline{\mu}}
\newcommand{\cX}{\mathbb{X}}

\begin{document}

\begin{center}
{\LARGE\bf Actions of monoidally equivalent compact quantum groups
and applications to probabilistic boundaries}

{\sc by An De Rijdt and  Nikolas Vander Vennet}
\end{center}

{\footnotesize  Department of Mathematics; K.U.Leuven;
Celestijnenlaan 200B; B--3001 Leuven (Belgium)

e-mail: an.derijdt@wis.kuleuven.be,
nikolas.vandervennet@wis.kuleuven.be}

\begin{abstract}
\noindent In \cite{BDRV}, Bichon, De Rijdt and Vaes introduced the
notion of monoidally equivalent compact quantum groups. In this
paper we prove that there is a natural bijective correspondence
between actions of monoidally equivalent quantum groups on unital
$C^*$-algebras or on von Neumann algebras. We apply this correspondence to study the behavior
of Poisson and Martin boundaries under monoidal equivalence of quantum groups.
\end{abstract}
\section*{Introduction}

After Woronowicz had introduced the notion of a compact quantum
group as a generalization of a compact group, many research topics
applying to compact groups were expanded to the general framework of
compact quantum groups. One of these topics concerns the study of
(ergodic) actions of compact groups on unital $C^*$-algebras (an
action on a $C^*$-algebra is ergodic if the fixed point algebra
reduces to the scalars).  We refer to the articles of H\o egh-Krohn,
Landstad  and St\o rmer \cite{HKLS} and Wasserman
\cite{wasser1,wasser2,wasser3} for a deep study of this topic. The
abstract theory of \emph{ergodic actions} of compact quantum groups
on \cst-algebras was initiated by Boca \cite{boca} and Landstad
\cite{landstad}. It turns out that the general theory of (ergodic)
actions of compact quantum groups on $C^*$-algebras is different
from the classical theory and in fact much richer. One major
difference is that the multiplicity of irreducible representations
in an ergodic action can be strictly greater than the
dimension of the representation space, which is impossible in the
classical case, where the dimension of the representation space is
actually an upper bound of this multiplicity. In the quantum case,
the upper bound is given by the quantum dimension which is usually
larger than the usual dimension.

In \cite{BDRV}, Bichon, the first author and Vaes introduced and
developed the notion of monoidally equivalent quantum groups. By
definition, compact quantum groups are called monoidally
equivalent if their representation categories are equivalent as
monoidal categories. In their article, they were able to describe
certain ergodic actions as unitary fiber functors on the
representation category. These ergodic actions are exactly the
ergodic actions of full quantum multiplicity. This provides us
with a powerful categorical tool for constructing  ergodic
actions. Moreover, these ergodic actions of full quantum
multiplicity provided the first examples of ergodic actions where
the multiplicity of the irreducible representations is strictly
greater than the dimension of the representation space.

In, \cite{Pin}, Pinzari and Roberts obtained a categorical
description of all ergodic actions of a compact quantum group.
Inspired by \cite{BDRV}, they describe an ergodic action (not
necessarily of full quantum multiplicity) of a compact quantum
group as a special kind of functor on the representation category.
When the ergodic action is of full quantum multiplicity, the
corresponding functor is just a unitary fiber functor as in
\cite{BDRV}. This categorical description yields a bijective
correspondence between ergodic actions of monoidally equivalent
quantum groups on unital $C^*$-algebras.

In this article, we obtain a bijective correspondence between (not
necessarily ergodic) actions of monoidally equivalent compact
quantum groups on unital $C^*$-algebras. Moreover, the
correspondence is of that kind that it preserves the spectral
subspaces of the actions. Restricting this to ergodic actions,
this just means that the multiplicities of the irreducible
representations are preserved through this correspondence. It
should be emphasized that our approach is not categorical. The
correspondence is obtained in a concrete, constructive way.

A major application of the bijective correspondence between
actions of monoidally equivalent quantum groups is found in the
study of Poisson and Martin boundaries for discrete quantum groups.
These boundaries find their origin in the study of random walks on discrete groups. A nice survey can be found in \cite{K1}. The study of random walks of discrete quantum groups was started by Biane, who considered duals of compact groups and obtained a theory which was parallel to the theory of random walks on discrete abelian groups.
Random walks on arbitrary discrete quantum groups, and their Poisson boundaries were
introduced by Izumi in \cite{iz1}, whose main motivation came from the study of infinite product actions of compact quantum groups. In \cite{iz1}, Izumi identified the
Poisson boundary of the dual of $\SU_q(2)$\ with the Podle\'s sphere \cite{Podles}. Later, Neshveyev and Tuset  \cite{NT}, associated a
Martin boundary to a random walk on a discrete quantum group and proved that the Martin boundary of the dual of $\SU_q(2)$\ is also given by the
Podle\'s sphere. In \cite{INT}, Izumi, Nesveyev and Tuset identified the Poisson boundary of $SU_q(n)$\ but its Martin boundary remains unknown.

Very recently, Tomatsu managed to identify the
Poisson boundaries of all amenable discrete quantum groups $\cGh$\
when its underlying compact quantum group $\cG$\ has commutative fusion rules. This Poisson boundary
appears to be the homogeneous space of $\cG$\ with respect to the
maximal closed quantum subgroup of Kac type \cite{tomatsu1}. This covers the results of Izumi, Neshveyev and Tuset, but is obtained in a completely different way and is also much more general. For Martin boundaries, there are no such results.

Considering the non-amenable case, Vaes and the second author \cite{VaeVa} have identified the Poisson and Martin boundary of the class of universal quantum groups $A_o(F)$\ with higher dimensional Podl\`{e}s spheres. When $\dim(F)\geq 3$, the dual of $A_o(F)$\ is not amenable.

We prove in this article a very general result. We provide a systematic
 method to relate Poisson and Martin boundaries for
the duals of monoidally equivalent quantum groups. The relation
goes as follows. Both the Poisson and Martin boundary of a discrete quantum group
$\cGh$, which is the dual of the compact quantum group $\cG$\
, admit a natural action of $\cG$. If the compact
quantum groups $\cG_1$ and $\cG_2$ are monoidally equivalent, the
boundaries of
 their duals $\cGh_1$ and $\cGh_2$ are related through the bijective
correspondence we obtained between the actions of $\cG_1$ and
$\cG_2$. This means that if we know the Poisson (Martin) boundary of the
dual of a compact quantum group $\cG$, we at the same time know it
for the duals of all compact quantum groups which are monoidally
equivalent with $\cG$.

Combining our result with Tomatsu's work, we give a concrete
identification of the Poisson boundary of a large class of discrete quantum groups. This method makes it also possible to obtain more examples of identifications of Poisson boundaries of
non-amenable discrete quantum groups. The main observation is that amenability is not preserved under monoidal equivalence. A first class
of examples of this kind are the universal orthogonal quantum
groups $A_o(F)$. If the dimension
 of $F$\ is greater than $3$, then $A_o(F)$\ is not coamenable. The
  quantum groups $A_o(F)$\ and $SU_q(2)$\ are monoidally
 equivalent for the right $q$. Moreover, the Poisson boundary of $SU_q(2)$\ was identified by Izumi and also,
 in a different way, by Tomatsu ($SU_q(2)$\ is coamenable). The correspondence just described gives a concrete identification of the Poisson boundary of
 $A_o(F)$. As we already saw, this result was already obtained by Vaes and the second author by another method \cite{VaeVa}.
  These were the first examples
 of identifications of Poisson boundaries of non-amenable discrete quantum groups.

 A second and new class of
 examples of the above type come from quantum automorphism
 groups $A_{aut}(D,\om)$, with $D$\ a finite dimensional $C^*$-algebra. These quantum groups have the fusion rules of $SO(3)$\ and are coamenable if and only if the dimension of the \cst-algebra
 is less than or equal to $4$. We prove, in that case that the maximal subgroup of Kac type is the one-dimensional torus $\mathbb{T}$. In combination with the result of Tomatsu, this provides us with an identification of its Poisson boundary. Using the fact that every quantum automorphism
 group is monoidally equivalent with a coamenable one, we obtain also an explicit
 identification of the Poisson boundary of the duals of all such quantum automorphism groups.

Because every $A_o(F)$\ is monoidally equivalent with an $SU_q(2)$, the correspondence of Martin boundaries under monoidally equivalent quantum groups gives a direct method to identify the Martin boundary of the duals of the universal compact quantum groups $A_o(F)$. This identification was already obtained by Vaes and the second author in \cite{VaeVa} by a different method using a result of \cite{Vbo}, allowing to deduce the Martin boundary, in the case of $A_o(F)$, from the Poisson boundary.

Finally, we would like to thank Stefaan Vaes for the numerous
remarks and careful reading of the manuscript.

\section{Notations}

Consider a subset $S$ of a \cst-algebra. We denote by $\langle S
\rangle$ the linear span of $S$ and by $[S]$ the closed linear span
of $S$. We use the notation $\om_{\eta,\xi}(a) = \langle \eta,a \xi
\rangle$ and we use inner products that are linear in the second
variable. Moreover we denote by $\xi^*:H\to\C:\eta\mapsto \langle
\xi,\eta\rangle$ and denote by $\overline{H}$ the dual Hilbert space
of $H$, i.e. $\overline{H}:=\{\xi^*\mid\xi\in H\}$.

The symbol $\ot$ denotes \emph{tensor products} of Hilbert spaces
and \emph{minimal} tensor products of \cst-algebras. We use the
symbol $\otalg$ for algebraic tensor products of $^*$-algebras and
$\otvn$ for the tensor product of von Neumann algebras. We also make
use of the leg numbering notation in multiple tensor products: if $a
\in A \ot A$, then $a_{12},a_{13},a_{23}$ denote the obvious
elements in $A \ot A \ot A$, e.g.\ $a_{12} = a \ot 1$.

The adjointable operators between \cst-modules or bounded
operators between Hilbert-spaces $H$ and $K$ are denoted by
$\cL(H,K)$. We also denote $\mathcal{L}(K,K)$\ by
$\mathcal{L}(K)$.

Let $\cB$ be a unital *-algebra. We call a linear map $\om:\cB\to\C$
such that $\om(1)=1$ a \emph{faithful state} if $\om(a^*a)\geq 0$
for all $a\in \cB$ and $\om(a^*a)=0$ if and only if $a=0$.

\section{Preliminaries}

\subsection*{Compact quantum groups}

We give a quick overview of the theory of compact quantum groups
which was developed by Woronowicz in \cite{wor2}. We refer to
\cite{MVD} for a survey of basic results.

\begin{definition}
A \emph{compact quantum group} $\cG$ is a pair $(C(\cG),\de)$, where
\begin{itemize}
\item $C(\cG)$ is a unital \cst-algebra; \item $\de : C(\cG) \recht C(\cG) \ot C(\cG)$ is a unital
  $^*$-homomorphism satisfying the \emph{co-associativity} relation
$$(\de \ot \id)\de = (\id \ot \de)\de \; ;$$
\item $\cG$ satisfies the \emph{left and right cancellation property}
  expressed by
$$\de(C(\cG))(1 \ot C(\cG)) \quad\text{and}\quad \de(C(\cG))(C(\cG) \ot 1)
\quad\text{are total in}\;\; C(\cG) \ot C(\cG) \; .$$
\end{itemize}
\end{definition}

\begin{remark}
The notation $C(\cG)$  suggests the analogy with the basic example
given by continuous functions on a compact group. In the quantum
case however, there is no underlying space $\cG$ and $C(\cG)$ is a
non-abelian \cst-algebra.
\end{remark}

A fundamental result in the theory of compact quantum groups is the
existence of a unique Haar state.

\begin{theorem}[Woronowicz, \cite{wor4}]
Let $\cG$ be a compact quantum group. There exists a unique state
$h$ on $C(\cG)$ which satisfies $(\id \ot h)\de(a) = h(a)1 = (h \ot
\id)\de(a)$ for all $a \in C(\cG)$. The state $h$ is called the
\emph{Haar state} of $\cG$.
\end{theorem}
Another crucial set of results in the framework of compact quantum
groups is the Peter-Weyl representation theory.
\begin{definition}
A \emph{unitary representation} $U$ of a compact quantum group $\cG$
on a Hilbert space $H$ is a unitary element $U \in \mathcal{L}(H \ot
C(\cG))$ satisfying
\begin{equation} \label{eq.rep}
(\id \ot \de)(U) = U_{12} U_{13} \; .
\end{equation}
Whenever $U^1$ and $U^2$ are unitary representations of $\cG$ on the
respective Hilbert spaces $H_1$ and $H_2$, we define
$$\Mor(U^1,U^2) := \{ T \in \cL(H_2,H_1) \mid U_1(T \ot 1) = (T \ot
1)U_2 \}\; .$$ The elements of $\mbox{Mor}(U^1,U^2)$ are called
\emph{intertwiners}. We use the notation $\mbox{End}(U) :=
\mbox{Mor}(U,U)$. A unitary representation $U$ is said to be
\emph{irreducible} if $\mbox{End}(U) = \C1$. If
$\mbox{Mor}(U^1,U^2)$
 contains a unitary
operator, the representations $U^1$ and $U^2$ are said to be
\emph{unitarily equivalent}.
\end{definition}

We have the following essential result.

\begin{theorem}
Every irreducible representation of a compact quantum group is
finite-dimensional. Every unitary representation is unitarily
equivalent to a direct sum of irreducibles.
\end{theorem}

Because of this theorem, we almost exclusively deal with
finite-dimensional representations. By choosing an orthonormal basis
of the Hilbert space $H$, a finite-dimensional unitary
representation of $\cG$ can be considered as a unitary matrix
$(U_{ij})$ with entries in $C(\cG)$ and \eqref{eq.rep} becomes
$$\de(U_{ij}) = \sum_k U_{ik} \ot U_{kj} \; .$$
The product in the \cst-algebra $C(\cG)$ yields a tensor product on
the level of unitary representations.
\begin{definition}
Let $U^1$ and $U^2$ be unitary representations of $\cG$ on the
respective Hilbert spaces $H_1$ and $H_2$. We define the tensor
product
$$U^1 \fancyot U^2 := U^1_{13} U^2_{23} \in \mathcal{L}(H_1 \ot H_2 \ot
C(\cG)) \; .$$
\end{definition}

\begin{notation}
Let $\cG$ be a compact quantum group. We denote by
$\mbox{Irred}(\cG)$ the set of equivalence classes of irreducible
unitary representations.  We choose representatives $U^x$ on the
Hilbert space $H_x$ for every $x \in \mbox{Irred}(\cG)$. Whenever
$x,y \in \mbox{Irred}(\cG)$, we use $x \ot y$ to denote the unitary
representation $U^x \fancyot U^y$. The class of the trivial unitary
representation is denoted by $\varepsilon$. We define the natural
numbers $\mbox{mult}(z,x \ot y)$ such that
$$x \ot y \cong \bigoplus_{z \in \mbox{Irred}(\cG)} \mbox{mult}(z,x \ot y) \cdot U^z \; .$$
The collection of natural numbers $\mult(z,x \ot y)$ are called the
\emph{fusion rules} of $\cG$.
\end{notation}

The set $\mbox{Irred}(\cG)$ is equipped with a natural involution $x
\mapsto \ox$ such that $U^{\overline{x}}$ is the unique (up to
unitary equivalence) irreducible unitary representation satisfying
$$\mbox{Mor}(x \ot \overline{x},\veps) \neq \{0\} \neq \mbox{Mor}(\overline{x} \ot x,\veps) \; .$$
The unitary representation $U^{\ox}$ is called the
\emph{contragredient} of $U^x$.

For every $x \in \mbox{Irred}(\cG)$, we take non-zero elements
${t_x}\in\mbox{Mor}(x\ot\ox,\veps)$ and $s_x
\in\mbox{Mor}(\overline{x}\ot x,\varepsilon)$ satisfying $(t_x^*\ot
1)(1\ot s_x)=1$. Write the antilinear map
\begin{equation}\label{eq.adjointmap}
j_x : H_x\to H_{\overline{x}}:\xi\mapsto (\xi^*\ot 1)t_x
\end{equation}
and define $Q_x := j_x^*j_x$\label{qqmatrix}. We normalize $t_x$
in such a way that $\Tr(Q_x)=\Tr(Q_x^{-1})$. This uniquely
determines $Q_x$ and fixes $t_x,s_x$ up to a number of modulus
$1$. Note that $t_x^*t_x = \Tr(Q_x)$.

\begin{definition}
For $x \in \mbox{Irred}(\cG)$, the value $\Tr(Q_x)$ is called the
\emph{quantum dimension} of $x$ and denoted by $\dimq(x)$. Note that
$\dimq(x) \geq \dim(x)$, with equality holding if and only if $Q_x =
1$.
\end{definition}

The irreducible representations of $\cG$ and the Haar state $h$ are
connected by the \emph{orthogonality relations}.
\begin{equation}\label{eq.orthogonality}
(\id \ot h)(U^x (\xi\eta^* \ot 1) (U^y)^*) =
\frac{\sde_{x,y}1}{\dimq(x)} \la \eta,Q_x\xi \ra \quad , \quad (\id
\ot h)((U^x)^*(\xi\eta^* \ot 1) U^y) = \frac{\sde_{x,y}1}{\dimq(x)}
\la \eta,Q^{-1}_x\xi \ra  \;
\end{equation}
for $\xi\in H_x$ and $\eta\in H_y$.

\begin{notation}
Let $\cG= (C(\cG),\de)$ be a compact quantum group. We denote by
$\mathcal{C}(\cG)$ the set of coefficients of finite dimensional
representations of $\cG$. Hence,
\[ \mathcal{C}(\cG)=\la (\om_{\xi,\eta}\ot\id)(U^x)\mid x\in\mbox{Irred}(\cG),\ \xi,\eta\in H_x
\ra \] Then, $\mathcal{C}(\cG)$ is a unital dense $^*$-subalgebra of
$C(\cG)$. Restricting $\de$ to $\mathcal{C}(\cG)$,
$\mathcal{C}(\cG)$ becomes a Hopf $^*$-algebra.\\
Also, for  $x\in \mbox{Irred}(\cG)$, denote by
$$\mathcal{C}(\cG)_x=\la (\om_{\xi,\eta}\ot\id)(U^x)\mid
\xi,\eta\in H_x \ra$$ Note that $\Delta: \mathcal{C}(\cG)_x\recht
\mathcal{C}(\cG)_x\ot \mathcal{C}(\cG)_x$\ and that
$\mathcal{C}(\cG)_x^*=\mathcal{C}(\cG)_{\overline{x}}$.
\end{notation}
\begin{definition} The \emph{reduced \cst-algebra} $C_r(\cG)$ is
  defined as the norm closure of $\mathcal{C}(\cG)$\
in the GNS-representation with respect to $h$. The \emph{universal
  \cst-algebra} $C_u(\cG)$ is defined as the enveloping \cst-algebra
of $\mathcal{C}(\cG)$. The \emph{von Neumann algebra}
$L^\infty(\cG)$ is defined as the von Neumann algebra generated by
$C_r(\cG)$. Note that if $\cG$ is the dual of a discrete group
$\Gamma$, we have $C_r(\cG) = C^*_r(\Gamma)$ and $C_u(\cG) =
C^*(\Gamma)$\ and $L^\infty(\cG)=\mathcal{L}(\Gamma)$.
\end{definition}
\begin{remark} \label{hopf} Given an arbitrary compact quantum group $\cG$, we have surjective homomorphisms
$C_u(\cG)\recht C(\cG)\recht C_r(\cG)$, but most of the time we
are only interested in $C_r(\cG)$ and $C_u(\cG)$. So, given the
underlying Hopf$^*$-algebra, there exists different
$C^*$-versions. From this point of view, we only consider two
quantum groups different if the underlying Hopf$^*$-algebras are
different.
\end{remark}

\begin{definition}\label{amena}
A compact quantum group $\cG$ is said to be \emph{coamenable}  if
the homomorphism $C_u(\cG) \recht C_r(\cG)$ is an isomorphism.
\end{definition}

\begin{proposition}\label{KMSh} The Haar state $h$ is a KMS-state on both $C_r(\cG)$\
  and $C_u(\cG)$ and the modular group is determined by
\begin{equation}\label{KMSh}
  (\id\ot
\sigma_t^h)(U^x)=(Q^{it}_x\ot 1)U^x(Q^{it}_x\ot 1)
\end{equation}
 for every $x\in
\Irred(\cG)$.
\end{proposition}

\subsection*{Discrete quantum groups and duality}

Following Van Daele (\cite{VD1}), a discrete quantum group is a
multiplier Hopf *-algebra whose underlying *-algebra is a direct
sum of matrix algebras. The dual of a compact quantum group is
such a discrete quantum group and is defined as follows.

\begin{definition}
Let $\cG$ be a compact quantum group. We define the dual (discrete)
quantum group $\cGh$ as follows.
\begin{equation*}
c_0(\cGh) = \bigoplus_{x \in \Irred(\cG)} \cL(H_x) \; , \qquad
\ell^\infty(\cGh) = \prod_{x \in \Irred(\cG)} \cL(H_x) \; .
\end{equation*}
We denote the minimal central projections of $\ell^\infty(\cGh)$ by
$p_x$, $x \in \Irred(\cG)$. We have a natural unitary $\Vb \in
\M(c_0(\cGh) \ot C(\cG))$ given by
\begin{equation}\label{UM}
\Vb = \bigoplus_{x \in \Irred(\cG)} U^x \; .
\end{equation}
This unitary $\Vb$ implements the duality between $\cG$ and $\cGh$.
We have a natural comultiplication
$$\deh : \ell^\infty(\cGh) \recht
\ell^\infty(\cGh) \otvn \ell^\infty(\cGh) : (\deh \ot \id)(\Vb) =
\Vb_{13} \Vb_{23} \; .$$
\end{definition}

One can deduce from this the following equivalent way to define the
coproduct structure on $\ell^\infty(\cGh)$.
$$\deh(a) S =
S a \quad\text{for all}\;\; a \in \ell^\infty(\cGh),\  S \in \Mor(y
\ot z,x) \; .$$  The notation introduced above is aimed to suggest
the basic example where $\cG$ is the dual of a discrete group
$\Gamma$, given by $C(\cG) = C^*(\Gamma)$ and $\de(\lambda_x) =
\lambda_x \ot \lambda_x$ for all $x \in \Gamma$. The map $x \mapsto
\lambda_x$ yields an identification of $\Gamma$ and $\Irred(\cG)$
and then, $\ell^\infty(\cGh) = \ell^\infty(\Gamma)$.

The discrete quantum group $\ell^\infty(\cGh)$ comes equipped with a
natural modular structure.

\begin{notation} \label{not.states}
 We have canonically defined states $\vphi_x$ and $\psi_x$ on $\cL(H_x)$
related to \eqref{eq.orthogonality} as follows.
\begin{align}
\psi_x(A)1 &= \frac{1}{\dimq(x)}t_x^* (A \ot 1) t_x = \frac{\Tr(Q_xA)}{\Tr(Q_x)}1 = (\id \ot h)(U^x (A \ot 1) (U^x)^*) \quad\text{and}\notag \\
\vphi_x(A)1 &=\frac{1}{\dimq(x)} t_{\ox}^* (1 \ot A) t_{\ox} =
\frac{\Tr(Q_x^{-1}A)}{\Tr(Q_x^{-1})}1 = (\id \ot h)((U^x)^*(A \ot 1)
U^x) \; ,\label{vallavalla}
\end{align}
for all $A \in \cL(H_x)$.
\end{notation}
\begin{remark}
The states $\vphi_x$ and $\psi_x$ are significant, since they
provide a formula for the invariant weights on $\ell^\infty(\cGh)$.
The left invariant weight is given by $\sum_{x \in \Irred(\cG)}
\dimq(x)^2 \psi_x $, and the right invariant weight is
 given by $\sum_{x \in \Irred(\cG)} \dimq(x)^2 \vphi_x$.
 \end{remark}

\begin{definition}A discrete quantum group $\cGh$\ is \emph{amenable} if there exists a left invariant mean on $\plon$,
 i.e. a state $m\in \plon^*$\ s.t.
$$m((\om\ot \id)\widehat{\de}(x))=m(x)\om(1)$$ for all $\om\in
\plon_*$\ and $x\in \plon$.
\end{definition}

\begin{remark} It was proven \cite{tomatsu2} that $\cGh$\ is \emph{amenable} if and only if $\cG$\ is \emph{coamenable}.
\end{remark}

\subsection*{Examples: the universal orthogonal compact quantum groups}

We consider a class of compact quantum groups which was introduced
by Wang and Van Daele in \cite{VDW}. These compact quantum groups
can in general not be obtained as deformations of classical objects.

\begin{definition}
Let $F \in \GL(n,\C)$ satisfying $F \overline{F} = \pm 1$. We define
the compact quantum group $\cG = A_o(F)$ as follows.
\begin{itemize}
\item $C_u(\cG)$ is the universal \cst-algebra with generators
  $(U_{ij})$ and relations making $U = (U_{ij})$ a unitary element of
  $\M_n(\C) \ot C(\cG)$ and $U = F \overline{U} F^{-1}$, where $(\overline{U})_{ij} = (U_{ij})^*$.
\item $\de(U_{ij}) = \sum_k U_{ik} \ot U_{kj}$.
\end{itemize}
\end{definition}

In these examples, the unitary matrix $U$ is a representation,
called the \emph{fundamental representation}. The definition of
$\cG=A_o(F)$ makes sense without the requirement $F \overline{F} =
\pm 1$, but the fundamental representation is irreducible if and
only if $F \overline{F} \in \R 1$. We then normalize such that $F
\overline{F} = \pm 1$.

\begin{remark}
It is easy to classify the quantum groups $A_o(F)$. For $F_1,F_2 \in
\GL(n,\C)$ with $F_i \overline{F}_i = \pm 1$, we write $F_1 \sim
F_2$ if there exists a unitary matrix $v$ such that $F_1 = v F_2
v^t$, where $v^t$ is the transpose of $v$. Then, $A_o(F_1) \cong
A_o(F_2)$ if and only if $F_1 \sim F_2$. It follows that the
$A_o(F)$ are classified up to isomorphism by $n$, the sign of $F
\overline{F}$ and the eigenvalue list of $F^*F$ (see e.g.\ Section 5
of \cite{BDRV} where an explicit fundamental domain for the relation
$\sim$ is described).

If $F \in \GL(2,\C)$, we get up to equivalence, the matrices
\begin{equation}\label{eq.fq}
F_q = \begin{pmatrix} 0 & |q|^{1/2} \\ - (\mbox{sgn} q) |q|^{-1/2} &
0
\end{pmatrix}
\end{equation}
for $q \in \, [-1,1]$, $q \neq 0$, with corresponding quantum
groups $A_o(F_q) \cong \SU_{q}(2)$, see \cite{wor}. In this case
the quantum dimension of the fundamental representation equals
$\Tr(F_q^*F_q)=\vert q+1/q\vert$.
\end{remark}

The following result has been proven by Banica \cite{banica2}. It
tells us that the compact quantum groups $A_o(F)$ have the same
fusion rules as the group $\SU(2)$.

\begin{theorem}
Let $F \in \GL(n,\C)$ and $F \overline{F} = \pm 1$. Let $\cG =
A_o(F)$. Then  $\Irred(\cG)$ can be identified with $\N$ in such a
way that
$$x \ot y \cong |x-y| \oplus (|x-y|+2) \oplus \cdots \oplus (x+y) \;
,$$ for all $x,y \in \N$.
\end{theorem}

Further on, we will introduce another class of compact quantum
groups that we need in this article, namely quantum
automorphism groups, but therefore we need the notion of an action.

\section{Actions of quantum groups}

\subsection*{Actions and spectral subspaces}
\begin{definition}\label{actie} Let $B$ be a unital \cst-algebra. A (right) action
  of $\cG$ on $B$ is a unital $^*$-homomorphism
$\sde : B \recht B \ot C(\cG)$ satisfying
\[(\sde\ot \id)\sde = (\id \ot \de)\sde \quad\text{and}\quad [\sde(B)(1\ot C(\cG))]=B\ot C(\cG)\; .\] The
action $\sde$ is said to be ergodic if the fixed point algebra
$B^{\sde}:=\{x\in B\mid  \sde (x)=x\ot 1\}$ equals $\C 1$. In that
case, $B$ admits a unique invariant state $\om$ given by $\om(b) 1 =
(\id \ot h)\sde(b)$.
\end{definition}

\begin{definition} \label{def.spectral}
Let $\delta : B \recht B \ot C(\cG)$ be an action of the compact
quantum group $\cG$ on the unital \cst-algebra $B$.
 For every $x \in \mbox{Irred}(\cG)$,  we define the spectral subspace associated with $x$ by
$$K_x = \{X \in \overline{H}_x \ot B \mid (\id \ot \sde)(X) = X_{12} U^x_{13} \} \;
.$$
\end{definition}
 Defining $\Hom(H_x,B)= \{ S:H_x \to B\mid S \ \textrm{linear and}\ \sde(S \xi )=(S \ot \id )(U^x(\xi \ot 1))\}$, we
have $K_x \cong \Hom(H_x,B)$, associating to every $X\in K_x$ the
operator $S_X:H_x\to B:\xi\mapsto X(\xi\ot 1)$.

\begin{remark}\label{Hilbert}
For each $x\in\mbox{Irred}(\cG)$, $K_x$ is  a bimodule over the
fixed point algebra $B^{\sde}$\ in a natural way. Indeed, for $a\in B^{\sde}$
and $X\in K_x$, $a\bullet X:=(1\ot a)X$ and $X\bullet a=X(1\ot a)$\
turns $K_x$ into a $B^{\sde}$-bimodule. Moreover, one can check
easily that
\begin{equation}\label{hil1}
\langle\cdot,\cdot \rangle:K_x \times K_x\to B^{\sde}:\langle
X,Y\rangle=XY^*\;.
\end{equation}
 gives an inner product, turning $K_x$
in a left Hilbert \cst-module over the fixed point algebra. We refer
to \cite{lance} for the theory of Hilbert C$^*$-modules.

We can also turn $K_x$ in a right Hilbert \cst-module. Denote by $E: B\to B^{\sde}: x\mapsto (\id\ot
h)\sde$ the conditional expectation onto the fixed point algebra.
For $X,Y\in K_x$, one can check, using the fact that
$(E\ot\id)\sde(x)=E(x)\ot 1$, that for each state $\om$ on
$B^{\sde}$, $(\id\ot\om E)(X^*Y)$ is an intertwiner for $U^x$ and
hence scalar. This means that we can define
\begin{equation}\label{hil2}
\langle\cdot,\cdot\rangle_{\sim}: K_x\times K_x\to B^{\sde}\quad
\text{by}\quad 1\ot\langle X,Y\rangle_{\sim}:=(\id\ot E)(X^*Y)\;,
\end{equation}
 which makes $K_x$ a right Hilbert \cst-module over $B^{\sde}$.

 In the case where $\sde$\ is ergodic with invariant state $\om$,
  $K_x$\ can be turned in a Hilbert space because
$B^\sde=\mathbb{C}$, with scalar product defined by $\langle
X,Y\rangle_l 1=YX^*$ and $\langle
X,Y\rangle_{r}1=(\id\ot\om)(X^*Y)$. Remark that we switched orders
in the first scalar product to have conjugate linearity in the first
variable.

\end{remark}

\begin{definition} \label{def.densesub}
We define $\cB$ as the subspace of $B$ generated by the spectral
subspaces, i.e.\
$$\cB:= \la X(\xi\ot 1)\mid x\in\mbox{Irred}(\cG),\ X\in K_x,\ \xi\in
H_x \ra \; .$$ Also, we define
$$\cB_x:= \la X(\xi\ot 1)\mid \ X\in K_x,\ \xi\in
H_x \ra \; .$$ Note that $\delta : \cB_x \recht \cB_x \otalg
\mathcal{C}(\cG)_x$ and that $\cB_x^* = \cB_{\overline{x}}$.
\end{definition}

Observe that $\cB$ is a dense unital *-subalgebra of $B$ and that
the restriction $\sde:\cB\to\cB\otalg \mathcal{C}(\cG)$ defines an
action of the Hopf $^*$-algebra $(\mathcal{C}(\cG),\de)$ on $\cB$.

\begin{remark} If $\delta$ is ergodic, $\cB_x$ is finite dimensional and its
  dimension is of the form $\dim H_x \cdot \mult(\delta,x)$, where
  $\mult(\delta,x)$ is called the multiplicity of $x$ in $\delta$.
  Note that as a vector space $B_x\simeq H_x\ot K_x$, so
  $\mult(\delta,x)=\dim K_x$.
\end{remark}

Suppose now that $\sde:B\recht B\ot C(\cG)$\ is an ergodic action.
Let $x \in \Irred(\cG)$. Take $t \in \Mor(\ox \ot x,\varepsilon)$,
normalized in such a way that $t^* t = \dimq(x)$. Define the
antilinear map
\begin{equation}\label{teees}
R_x : K_x \recht K_{\ox} : R_x(v) = (t^* \ot 1)(1 \ot v^*) \; .
\end{equation}
Since $t$ is fixed up to a number of modulus one, $L_x := R_x^* R_x$
is a well defined positive element of $\cL(K_x)$.

\begin{definition}\label{quam}
We put $\multq(x) := \sqrt{\Tr(L_x) \Tr(L_{\ox})}$ and we call
$\multq(x)$ the \emph{quantum
  multiplicity of $x$ in $\sde$}.
\end{definition}

\begin{remark}
It can be proven, for example in \cite{BDRV}, that $\multq(x) \leq
\dimq(x)$ for all $x \in \Irred(\cG)$. If equality holds for all
$x \in \cGh$, we say that $\sde$ is of \emph{full quantum
multiplicity}.
\end{remark}

\begin{terminology}
An action $\sde:B\to B\ot C(\cG)$ of $\cG$ on $B$ is said to be
\emph{universal} if $B$ is the universal enveloping \cst-algebra of
$\cB$. It is said to be \emph{reduced} if the conditional
expectation $(\id\ot h)\sde$ of $B$ on the fixed point algebra
$B^\sde$ is faithful.
\end{terminology}

\begin{remark} \label{prop.densesub} From remark \ref{hopf}, we saw that a compact quantum group $(C(\cG),\Delta)$\ has many
$C^*$-versions, while the underlying Hopf$^*$-algebra is the same.
The same remark applies to actions. We have that
$B_u\recht B\recht B_r$\ for an action $\sde:B\recht B\ot
C(\cG)$. So again, we only consider two actions to be different if
the underlying Hopf $^*$-algebra actions are different. We make
extensively use of this fact.
\end{remark}

Actions on von Neumann algebras are defined as follows.

\begin{definition} A right action of a compact (resp.\ discrete)
  quantum group $\cG$ (resp. $\widehat{\cG}$) on a von Neumann algebra $N$ is an injective
 normal unital $^*$-homomorphism
$$\sde:N\recht N\otvn L^{\infty}(\cG)\qquad \text{resp.}\quad \sde:N\recht N\otvn \ell^{\infty}(\widehat{\cG})$$
satisfying $(\sde\ot \id)\sde = (\id \ot \de)\sde $, resp.\ $(\delta
\ot \id)\delta = (\id \ot \deh)\delta$.
\end{definition}

\begin{remark}
In the case of an action of a compact
  quantum group on a von Neumann algebra, we do not require the density condition like for \cst-algebraic
actions. The reason is that this is automatically fulfilled for
von Neumann algebras. This is a quite deep result and we refer to
\cite{vaes}, theorem 2.6 for a proof. This implies that the
spectral subalgebra as defined in \ref{def.densesub} remains
(weakly) dense in $N$.
\end{remark}

\begin{remark} Because every action $\sde:B\recht B\ot C(\cG)$\ has an unitary implementation, it can be extended to a von Neumann algebraic action.
\end{remark}

\subsection*{Quantum subgroups and homogeneous spaces}

\begin{definition}\label{def.quantumsubgroup}
Let $(\cG,\de_\cG)$ and $(\cH,\de_\cH)$ be compact quantum groups.
We call $\cH$ a closed \emph{quantum subgroup} of $\cG$ whenever
there is given a surjective *-homomorphism $r_{\cH}:\cC(\cG)\to
\cC(\cH)$ satisfying $\de_{\cH}\circ r_{\cH}=(r_{\cH}\ot
r_{\cH})\de_{\cG}$.
\end{definition}

\begin{definition}\label{homogeen}
Let  $(\cG,\de_{\cG})$\ a compact quantum group with quantum
subgroup $(\cH,\de_{\cH})$. Define the Hopf*-algebra action
$\gamma_{\cH}:\cC(\cG)\to \cC(\cH)\otalg \cC(\cG):x\mapsto
(r_{\cH}\ot\id)\de_{\cG}(x)$. Define the \emph{homogeneous space}
$\cC(\cH\backslash\cG)$\ as the fixed point subalgebra of $\cC(\cG)$
under $\gamma_{\cH}$.
\end{definition}

\begin{remark}\label{qawz} The restriction of the comultiplication to
$\cC(\cH\backslash\cG)$\ gives a Hopf$^*$-action
\begin{equation}\label{lopi}
\de_{\cH\backslash\cG}:\cC(\cH\backslash\cG)\to
\cC(\cH\backslash\cG)\otalg \cC(\cG)\;.
\end{equation}
\end{remark}

Since the action $\gamma_\cH$ is invariant under de Haar measure of
$\cG$, we can extend it to $C_r(\cG)$ and $L^\infty(\cG)$ and hence
define $C_r(\cH\backslash\cG)$ and $L^\infty(\cH\backslash\cG)$. By
universality, $\gamma_\cH$ is also extendable to $C_u(\cG)$, which
gives us $C_u(\cH\backslash\cG)$.

\label{sqsqsq} The restriction of the comultiplication to
$C_r(\cH\backslash\cG)$, respectively $C_u(\cH\backslash\cG)$, or
$L^{\infty}(\cH\backslash\cG)$\ gives again an action as in
formula $(\ref{lopi})$.

\begin{lemma}\label{coeenheid}
The restriction of $r_{\cH}$ to the quotient $\cC(\cH\backslash\cG)$
is the co-unit $\eps_{\cG}$.
\end{lemma}
\begin{proof}
For $a\in \cC(\cH\backslash\cG)$,
\[\de_{\cH}(r_{\cH}(a))=(r_{\cH}\ot r_{\cH})\de_{\cG}(a)=1\ot r_{\cH}(a)\;,\]
 We now apply $(\id\ot\eps_{\cH})$ to
both sides of the equation and use the fact that
$\eps_{\cH}r_{\cH}=\eps_{\cG}$ (\cite{tomatsu1}). Then
\[r_{\cH}(a)=(\id\ot\eps_{\cH})\de_{\cH}(r_\cH(a))=\eps_{\cG}(a)1\;,\] which
ends the proof.
\end{proof}

In the last chapter, we will need a special kind of subgroup.

\begin{definition}\label{Kaccc}
Consider a compact quantum group $(\cG,\de_{\cG})$. We call a
quantum subgroup $(\cH,\de_{\cH})$ of Kac type maximal, if for any
quantum subgroup $\mathbb{K}$ of Kac type,
$L^{\infty}(\cH\backslash\cG)\subset
L^{\infty}(\mathbb{K}\backslash\cG)$.
\end{definition}

Every compact quantum group has a unique maximal quantum subgroup of
Kac type (see \cite{Soltan}). We call it the \emph{canonical Kac
subgroup} of the quantum group.

\subsection*{Invariant subalgebras}

A more general notion is that of an invariant subalgebra.

\begin{definition}
Consider a compact quantum group $\cG$ with comultiplication $\de$.
A \emph{right invariant subalgebra} of $\cG$ is a unital
\cst-algebra $B\subset C(\cG)$ such that $\de(B)\subset B\ot
C(\cG)$.
\end{definition}

We can define an ergodic action $\sde$ of $\cG$ on $B$ by just
restricting $\de$ to $B$. We get the following easy proposition.

\begin{proposition}\label{prop.coideal}
 Consider a compact quantum group $\cG$ and a right invariant subalgebra $B$ of $C(\cG)$. Denote
 the action of $\cG$ on $B$ by $\sde$. For all $x\in \Irred(\cG)$,
 $\mult(\sde,x)\leq\dim(x)$ and equality in all $x$ is only reached
 when $B=C(\cG)$.
\end{proposition}

\begin{proof}
Let $x\in Irred(\cG)$. From the definition of a spectral subspace,
we get
\[K_x=\{X \in \overline{H}_x \ot B \mid (\id \ot \de)(X) =
X_{12} U^x_{13} \}\;.\]It is clear that
\[K_x\subset\mathbb{K}_x:=\{X
\in \overline{H}_x \ot C(\cG) \mid (\id \ot \de)(X) = X_{12}
U^x_{13} \}\] with $\mathbb{K}_x$\ the spectral subspace of the
comultiplication $\de$. Now $\mathbb{K}_x\cong \oH_x$ where the
bijection is given by $\oH_x\to K_x:\overline{\xi}\mapsto (\xi^*\ot
1)U^x$. Then $\mult(\sde,x)=\dim(K_x)\leq\dim(H_x)=\dim(x)$.

Equality for all $x\in \Irred(\cG)$ means that $K_x=\mathbb{K}_x$,
so $B_x=\cC(\cG)_x$ and hence $B=C(\cG)$.
\end{proof}

\subsection*{Quantum automorphism groups}\label{qaut}

In this section we consider a class of universal quantum groups,
namely the quantum automorphism groups as studied by Wang in
\cite{wang2} and Banica in \cite{banica5,banica6}. We only
consider \cst-algebras with a special kind of states.

\begin{definition} Let $(D,\om)$ be a
finite dimensional \cst-algebra of dimension $\geq 4$ with a
state. Denote by $\mu:D\ot D\to D$ the
multiplication. Take $\sde>0$. If for the inner product implemented by
$\om$, $\mu\mu^*=\sde^2 1$, we call $\om$ a $\sde$-form.
\end{definition}

If $D$ is a matrix-algebra, every state is of the form $\Tr(F\cdot)$
and a $\sde$-form with $\sde^2=\Tr(F^{-1})$. This can easily be
checked by writing out $\mu\mu^*$ in terms of the orthonormal basis
$(e_{ij}F^{-\frac{1}{2}})_{i,j=1\cdots n}$ of $D$.

We can now give the definition of a quantum automorphism group:

\begin{definition}\cite{banica6}
Let $(D,\om)$ be an finite-dimensional \cst-algebra with a
$\sde$-form. We define the compact quantum group
$\cG=A_{aut}(D,\om)$ as follows. $\cG$ is defined by an action $\alpha:D\to D\ot
C_u(\cG)$ with the following properties:
\begin{itemize}
\item
$C_u(\cG)$ is defined as the universal \cst-algebra generated by
\[\{(\om\ot\id)\alpha(a)\mid \om\in D^*,\ a\in D\}\;.\]
\item Whenever $\beta:D\to D\ot C_u(\cG_1)$ is an action of a compact quantum group
$\cG_1$, there exists a unique *-homomorphism $\pi:C_u(\cG)\to C_u(\cG_1)$ satisfying
$\beta=(\id\ot\pi)\alpha$.
\end{itemize}
\end{definition}

\begin{remark} In this article, we consider only the cases
where $n\geq 4$. In the cases $n=1,2,3$, we just get the
permutation group $S_n$.
\end{remark}

\subsubsection*{Representation Theory}\label{banica}\label{structure}

In \cite{banica5}, Banica has determined the irreducible
representations and their fusion rules for all quantum automorphism
groups.

If $B$ and $\om$ are as above, the fusion rules of $A_{aut}(D,\om)$
are those of $\SO(3)$. This means that the irreducible
representations are labeled by $\N$. We choose $U_i\in \cL(H_i)\ot
\cC(A_{aut}(D,\om))$ the representative of the irreducible
representation with label $i$ in such a way that $U_0$ is the
trivial representation $\veps$ and that $U=U_0\oplus U_1\in\cL(D)\ot
\cC(A_{aut}(D,\om))$ is the fundamental representation. The fusion
rules are given by:
\[U_i\ot U_j=U_{|i-j|}+ U_{|i-j|+1}+\cdots + U_{i+j}\;.\]

\section {Monoidal equivalence}
\subsection*{General theory}
The notion of monoidal equivalence was introduced in  \cite{BDRV}.
In this section, we give an overview of the results we will need.
\begin{definition}[Def.\ 3.1 in \cite{BDRV}] \label{def.moneq}
Two compact quantum groups $\cG_1=(C(\cG_1),\de_1)$ and
$\cG_2=(C(\cG_2),\de_2)$ are said to be \emph{monoidally
equivalent} if there exists a bijection
$\vphi:\mbox{Irred}(\cG_1)\to\mbox{Irred}(\cG_2)$ satisfying
$\vphi(\varepsilon) = \varepsilon$, together with linear
isomorphisms
\[\vphi:\Mor(x_1 \ot \cdots \ot x_r ,y_1\ot\cdots\ot
y_k)\to\Mor(\vphi(x_1) \ot \cdots \ot \vphi(x_r),\vphi(y_1)\ot\cdots
\ot \vphi(y_k))\] satisfying the following conditions:
\begin{equation*}
\begin{alignedat}{2}
\vphi(1) &= 1 & \qquad
\vphi(S \ot T) &= \vphi(S) \ot \vphi(T) \\
\vphi(S^*) &= \vphi(S)^* & \qquad \vphi(S T) &=\vphi(S) \vphi(T)
\end{alignedat}
\end{equation*}
whenever the formulas make sense. In the first formula, we consider
$1 \in \Mor(x,x) = \Mor( x \ot \varepsilon,x) = \Mor(\varepsilon \ot
x,x)$. Such a collection of maps $\vphi$ is called a \emph{monoidal
equivalence} between $\cG_1$ and $\cG_2$.
\end{definition}

By Theorem 3.9 and Proposition 3.13  of \cite{BDRV}, we have the
following fundamental result.

\begin{theorem}\label{construction}
Let $\vphi$ be a monoidal equivalence between compact quantum groups
$\cG_1$ and $\cG_2$.
\begin{itemize}
\item There exist a unique unital $^*$-algebra $\cB$ equipped with a faithful state $\om$ and unitary elements $X^x \in
\B(H_x,H_{\vphi(x)}) \ot \cB$ for all $x \in \mbox{Irred}(\cG_1)$,
satisfying
\begin{enumerate}
\item $X^y_{13} X^z_{23} (S \ot 1) = (\vphi(S) \ot 1)X^x
  \quad\text{for all}\quad S \in \Mor(y\ot z,x) \; ,$
\item the matrix coefficients of the $X^x$ form a linear basis of
  $\cB$,
\item $(\id \ot \om)(X^x) = 0 \quad\text{if}\quad x \neq \veps$.
\end{enumerate}
\item There exists  unique commuting ergodic actions $\sde_1 : \cB \recht \cB \otalg \cC(\cG_1)$ and $\sde_2 : \cB \recht
\cC(\cG_2)\otalg \cB$ satisfying
$$(\id \ot \sde_1)(X^x) = X^x_{12} U^x_{13}\qquad \mbox{and}\qquad (\id \ot \sde_2)(X^x) = U^{\varphi(x)}_{12} X^x_{13}$$
for all $x \in \Irred(\cG)$. \item The state $\om$ is invariant
under $\sde_1$ and $\sde_2$. Denoting
  by $B_r$ the \cst-algebra generated by $\cB$ in the
  GNS-representation associated with $\om$ and denoting by $B_u$ the
  universal enveloping \cst-algebra of $\cB$, the actions $\sde_1,\sde_2$
  admit unique extensions to actions on $B_r$ and $B_u$.
\end{itemize}
\end{theorem}
This algebra $\cB$\ is called the link algebra of $\cG_1$\ and
$\cG_2$\ under the monoidal equivalence $\varphi$.

Note that in the case $\cG=\cG_1=\cG_2$ and $\vphi$ the identity
map, we have $\cB=\mathcal{C}(\cG)$ and $X^x=U^x$ for every $x\in
\Irred(\cG)$. The following unitary operator generalizes \eqref{UM}.
\begin{equation}\label{eq.X}
\mathbb{X}:=\bigoplus_{x\in\Irred(\cG)}X^x \quad\text{where}\quad
\mathbb{X} \in \prod_{x \in \Irred(\cG)}
\bigl(\cL(H_x,H_{\vphi(x)}) \ot B \bigr) \; .
\end{equation}

\begin{proposition}\label{KMSw}
The invariant state $\om$ is a KMS state on $B_r$ and $B_u$ and its
modular group is determined by
\begin{equation}\label{KMSS}
(\id\ot \sigma^{\om}_t)(X^x)=(Q_{\vphi(x)}^{it}\ot 1)X^x(Q_x^{it}\ot
1)
\end{equation}
for every $x\in \Irred(\cG_1)$.
\end{proposition}

\begin{remark}
Define $\cB_x:=\la (\om_{\xi,\eta}\ot \id)(X^x)\mid\ \xi\in
H_{\varphi(x)}, \eta\in H_x \ra \; .$ Then, as a vector space
$$\mathcal{B}=\bigoplus_{x\in \mbox{Irred}(\cG)}\cB_x\; .$$ Moreover, the $\cB_x$ are exactly the
spaces $\cB_x$\ in definition \ref{def.densesub}  coming from the
spectral subspaces of $\sde_1$ and $\sde_2$, while $\mathcal{B}$ is
exactly the dense $^*$-algebra given in Definition
\ref{def.densesub}.
\end{remark}

The orthogonality relations \eqref{eq.orthogonality} generalize and
take the following form.
\begin{equation}\label{eq.orthogonalitymon}
\begin{split}
(\id \ot \om)(X^x (\xi_1\eta_1^* \ot 1) (X^y)^*) &=
\frac{\sde_{x,y}1}{\dimq(x)} \la \eta_1,Q_x\xi_1 \ra \; , \\ (\id
\ot \om)((X^x)^*((\xi_2\eta_2^* \ot 1) X^y) &=
\frac{\sde_{x,y}1}{\dimq(x)}\la \eta_2,Q_{\varphi(x)}^{-1}\xi_2 \ra
\; ,
\end{split}
\end{equation}
for $\xi_1\in H_x$, $\eta_1\in H_y$, $\xi_2\in H_{\varphi(x)}$ and
$\eta_2\in H_{\varphi(y)}$.

\subsection*{Concrete examples}

In this section, we investigate in a closer way monoidal
equivalence for specific quantum groups, namely the universal
quantum groups $A_o(F)$\ and the quantum automorphism groups. The case of the quantum groups $A_o(F)$\ was already studied in detail in
\cite{BDRV}. If  $\cG_1=A_o(F_1)$ and
$\cG_2=A_o(F_2)$, the following theorem gives a concrete expression of their link algebra.

\begin{theorem}[Thms.\ 5.3 and 5.4 in \cite{BDRV}] \label{constructionk}
 Let $F_1\in M_{n_1}(\mathbb{C})$ and $F_2\in M_{n_2}(\mathbb{C})$ such that $F_1\overline{F}_1=\pm 1$\
 and $F_2\overline{F}_2=\pm 1$.
\begin{itemize}
\item  The compact quantum groups $A_o(F_1)$ and $A_o(F_2)$
  are monoidally equivalent iff $F_1\overline{F}_1$ and $F_2\overline{F}_2$ have the same sign and
  $\Tr(F_1^*F_1)=\Tr(F_2^*F_2)$.
\item Assume that $A_o(F_1)$ and $A_o(F_2)$
  are monoidally equivalent.
  Denote by $C_u(A_o(F_1,F_2))$ the universal unital \cst-algebra
  generated by the coefficients of
$$Y \in M_{n_2,n_1}(\C) \ot C_u(A_o(F_1,F_2)) \quad\text{with relations}\quad
Y \;\;\text{unitary}\quad\text{and}\quad Y = (F_2 \ot
1)\overline{Y} (F_1^{-1} \ot 1) \; .$$ Then, $C_u(A_o(F_1,F_2))
\neq 0$ and there exists a unique pair of commuting universal
ergodic actions, $\sde_1$ of $A_o(F_1)$ and $\sde_2$ of
$A_o(F_2)$, such that $$(\id \ot \sde_1)(Y) = Y_{12} (U_1)_{13}
\quad\text{and}\quad (\id \ot \sde_2)(Y) = (U_2)_{12} Y_{13} \;
.$$ Here, $U_i$ denotes the fundamental representation of
$A_o(F_i)$.
\item $(C_u(A_o(F_1,F_2)),\sde_1,\sde_2)$ is isomorphic with the \cst-algebra $B_u$ and the actions thereon given by
theorem \ref{construction}.
\end{itemize}
\end{theorem}

\begin{remark} It is also true that any compact quantum group
which is monoidally equivalent with $A_o(F)$\ where $F\in
\GL(n,\mathbb{C})$\ and $F\overline{F}=\pm 1$\ is itself of the
form of $A_o(F_1)$\ where $F_1\in \GL(n_1,\mathbb{C})$\ and
$F_1\overline{F_1}=\pm 1$. Even more holds, Banica \cite{banica2} showed that any quantum group with
fusion rules of $SU(2)$\ is of the form $A_o(F)$\  where $F\in
\GL(n,\mathbb{C})$\ and $F\overline{F}=\pm 1$.
\end{remark}

Next, we obtain a concrete expression of the link algebra in the
case that $\cG_1=A_{aut}(D_1,\om_1)$ and
$\cG_2=A_{aut}(D_2,\om_2)$. We prove the following theorem.

\begin{theorem}\label{extrabew}
Let $D_1$ and $D_2$ be finite dimensional \cst-algebras and
$\om_1$ and $\om_2$ respectively a $\sde_1$-form and a
$\sde_2$-form on $D_1$, respectively $D_2$.

\begin{itemize}
\item The compact quantum groups $\cG_1=A_{aut}(D_1,\om_1)$ and
$\cG_2=A_{aut}(D_2,\om_2)$ are monoidally equivalent if and only
if $\sde_1=\sde_2$. \item Suppose that $A_{aut}(D_1,\om_1)$ and
$A_{aut}(D_2,\om_2)$ are monoidally equivalent. Denote by
$C_u(A_{aut}((D_1,\om_1),(D_2,\om_2)))$ the universal \cst-algebra
gene\-rated by the matrix elements of a unital *-homomorphism
\[\gamma:D_1\to D_2\ot C_u(A_{aut}((D_1,\om_1),(D_2,\om_2)))\]
\[\text{ with relations } (\om_2\ot\id)\gamma(x)=\om_1(x)1\text{
for all }x\in D_1\;.\] Then
$C_u(A_{aut}((D_1,\om_1),(D_2,\om_2)))\neq 0$ and there exists a
unique pair of commuting ergodic actions of full quantum
multiplicity $\sde_1$ of $A_{aut}(D_1,\om_1)$ and $\sde_2$ of
$A_{aut}(D_2,\om_2)$, such that
\[(\id\ot\sde_1)\gamma=(\gamma\ot\id)\beta_1 \text{ and }
(\id\ot\sde_2)\gamma=(\beta_2\ot\id)\gamma\;,\] where
$\beta_1:D_1\to D_1\ot C_u(A_{aut}(D_1,\om_1))$ and
$\beta_2:D_2\to D_2\ot C_u(A_{aut}(D_2,\om_2))$ are the actions of
the quantum automorphism groups. \item
$(C_u(A_{aut}((D_1,\om_1),(D_2,\om_2))),\sde_1,\sde_2)$ is isomorphic with the
\cst-algebra $B_u$ and the actions thereon given by proposition
\ref{construction}.
\end{itemize}
\end{theorem}

\begin{proof}
Denote by $\mu_1$, $\mu_2$ and $\eta_1$, $\eta_2$ the
multiplication and unital map of respectively $D_1$ and $D_2$. The
proof of the first point goes as follows. First suppose that
$\sde_1=\sde_2$. Take now $U$, respectively $V$ the fundamental
representation of $A_{aut}(D_1,\om_1)$, respectively
$A_{aut}(D_2,\om_2)$ corresponding to the actions of this quantum
groups. Consider the graded \cst-algebras $(\Mor(U^m,U^n))_{n,m}$
and $(\Mor(V^m,V^n))_{n,m}$. We know from \cite{banica5} that
there is an isomorphism $\pi: (\Mor(U^n,U^m))_{n,m}\to
(\Mor(V^n,V^m))_{n,m}$ which satisfies $\pi(\mu_1)=\mu_2$ and
$\pi(\eta_1)=\eta_2$ We now can work analogously to the case of
$A_o(F)$ that was covered in \cite{BDRV}.

We now set $\Irred(\cG_1)=\N$ and $P_n\in\Mor(U^n,U^n)$ the unique
projection for which $P_nT=0$ for all $r<n$ and all
$T\in\Mor(U^r,U^n)$. We define $U_n$ as the restriction of $U^n$
to the image of $P_n$ and identify
\[\Mor(n_1\ot\cdots\ot n_r,m_1\ot\cdots\ot
m_k)\] \[=(P_{m_1}\ot\cdots\ot
P_{m_k})\Mor(U^{n_1+\cdots+n_r},U^{m_1+\cdots+m_k})(P_{n_1}\ot\cdots\ot
P_{n_r})\;.\] Define now $H_{\psi(n)}:=\pi(P_n)D_1^n$  and define
for $S\in\Mor(n_1\ot\cdots\ot n_r,m_1\ot\cdots\ot m_k)$, $\psi(S)$
by the restriction of $\pi$ to $\Mor(n_1\ot\cdots\ot
n_r,m_1\ot\cdots\ot m_k)$. Then $\psi$ is a unitary fiber functor
which gives a monoidal equivalence between $\cG_1$ and $\cG_2$.

Conversely, suppose that $A_{aut}(D_1,\om_1)\meq
A_{aut}(D_2,\om_2)$. Denote by $u_1$, respectively $v_1$ the
irreducible representation with label $1$ of $A_{aut}(D_1,\om_1)$
and $A_{aut}(D_2,\om_2)$. Then $\dim_q(u_1)1=\om_1
\mu_1\mu_1^*\om_1^*=\sde_1^2 1$ and because monoidal equivalence
preserves the quantum dimension, $\sde_1$ and $\sde_2$ must be
equal. This proofs the first part of the theorem.

For the proof of the other parts of the theorem, we first make the
following observation.  Consider two finite dimensional
\cst-algebras $(D_1,\om_1)$ and $(D_2,\om_2)$ with $\sde$-forms
and their quantum automorphism groups\\
$A_{aut}(D_1,\om_1):=\cG_1$ and $A_{aut}(D_2,\om_2):=\cG_2$.
Denote now by $H_1^i=D_i\ominus \C$ and $U_i^1\in
\cL(H_1^i)\otimes \cC(\cG_i)$ for $i=1,2$ the representative of
the irreducible representation with label $1$. Denote by
$\theta_i\in\Mor((U_i^1),(U_i^1)^2)$\ and
$\gamma_i\in\Mor(U_i^0,(U_i^1)^2)$\ the obvious "components" of
the multiplication. From the construction in the first part of the
theorem, it follows that there is a monoidal equivalence $\vphi$
between $\cG_1$ and $\cG_2$ which sends $\theta_1$\ and
$\gamma_1$, to $\theta_2$\ and $\gamma_2$. If we further below
talk about the monoidal equivalence between $A_{aut}(D_1,\om_1)$
and $A_{aut}(D_2,\om_2)$, we will always mean this one.

We first remark that if $C_u(A_{aut}((D_1,\om_1),(D_2,\om_2)))\neq
0$, the actions are given by universality. Indeed,
\[(\gamma\ot\id)\beta_1:D_1\to D_2\ot A_{aut}((D_1,\om_1),(D_2,\om_2))\ot
A_{aut}(D_1,\om_1)\] is a *-homomorphism which satisfies
\[(\om_2\ot\id\ot\id)(\gamma\ot\id)\beta_1(x)=\om_1(x)1\] for $x\in
D_1$. So by universality, there exists a *-homomorphism
\[\sde_1:A_{aut}((D_1,\om_1),(D_2,\om_2))\to
A_{aut}((D_1,\om_1),(D_2,\om_2))\ot A_{aut}(D_1,\om_1)\]
satisfying $(\id\ot\sde_1)\gamma=(\gamma\ot\id)\beta_1$. Because
$\beta_1$ is an action and the coefficients of $\gamma$ generate
$C_u(A_{aut}((D_1,\om_1),(D_2,\om_2))$, it follows that $\sde_1$
is an action. We define $\sde_2$ in an analogous way.

Consider now the \cst-algebra $B_u$ we get from the monoidal
equivalence. Denote by $\theta_i$, $\gamma_i$ the components of
the multiplication of $D_i$, $i=1,2$. As we said above, we may
suppose that the monoidal equivalence sends $\theta_1$ and
$\gamma_1$ respectively to $\theta_2$ and $\gamma_2$. Denote by
$U_i$ the irreducible representation of $A_{aut}(D_i,\om_i)$ with
label $1$. Because every irreducible representation is contained
in a tensor power of the one with label $1$, the matrix
coefficients of $X^1\in \cL(D_1\ominus\C,D_2\ominus\C)\ot B_u$
generate $B_u$ as a \cst-algebra. By identification, $X^1$
provides us with a linear map
\[\Gamma: D_1\ominus\C\to (D_2\ominus\C)\ot B_u\] which we can easily extend to $D_1$ by
setting $\Gamma(1)=1$. Because \[X^1(\theta_1\ot 1)=(\theta_2\ot
1)X^1_{13}X^1_{23} \text{ and }(\gamma_1\ot 1)=(\gamma_2\ot
1)X^1_{13}X^1_{23}\;,\] $\Gamma$ is multiplicative, obviously
unital and $\om_1(x)1=(\om_2\ot\id)\Gamma(x)$. It also preserves
the involution because $X_{23}(\gamma_1^*\ot
1)=X^*_{13}(\gamma_2^*\ot 1)$ and $\gamma_1$ and $\gamma_2$
implement the involution on respectively $D_1$ and $D_2$. By
universality there exists now a unital *-homomorphism
\[\rho:C_u(A_{aut}((D_1,\om_1),(D_2,\om_2)))\to B_u\] such that
$\Gamma=(\id\ot\rho)\gamma$. It is now left to show that $\rho$ is
an isomorphism.

Because $\gamma$ satisfies the equation
$(\om_2\ot\id)\gamma(x)=\om_1(x)1$, we can look at the restriction
of $\gamma$ given by
\[\gamma: D_1\ominus\C\to (D_2\ominus\C)\ot C_u(A_{aut}((D_1,\om_1),(D_2,\om_2)))\;.\]
Denote by $Y\in \cL((D_1\ominus \C), (D_2\ominus\C))\ot
C_u(A_{aut}((D_1,\om_1),(D_2,\om_2)))$ the element corresponding
to this restricted *-homomorphism. This element satisfies the
equations
\begin{align*}
Y(\theta_1\ot 1)&=(\theta_2\ot 1)Y_{13}Y_{23} \text{ and }\\
\gamma_1\ot 1&=(\gamma_2\ot 1)Y_{13}Y_{23}
\end{align*}
because $\gamma$ is a unital homomorphism. Remark that $Y$ is
unitary because $\gamma$ also preserves the involution. Because
the multiplication and the unital map generate all the
intertwiners of $A_{aut}(D_i,\om_i)$, $i=1,2$, and so
is also true for $\theta_i$ and $\gamma_i$, it holds that
\[Y^{\ot n}(P_n\ot 1)=(Q_n\ot 1)Y^{\ot n}\] where $P_n$ and $Q_n$ are the unique
projections in respectively $\Mor(U^n_1,U_1^n)$ and
$\Mor(U_2^n,U_2^n)$ on the irreducible representation with label
$n$. Defining $\sigma$ such that $(\id\ot\sigma)(X^n)=Y^{\ot
n}(P_n\ot 1)$, gives a unital *-homomorphism with
$\sigma\rho=\rho\sigma=id$.
\end{proof}

\begin{remark}
We can also prove that every compact quantum group $\cG$ which is
monoidally equivalent to a quantum automorphism group
$A_{aut}(D,\om)$ is isomorphic to another quantum automorphism
group $A_{aut}(D_1,\om_1)$. It is to our best knowledge not clear if
every compact quantum group with the fusion rules of $\SO(3)$ is a
quantum automorphism group.
\end{remark}

\section{The Poisson boundary of a discrete quantum group}

We give a brief survey of Izumi's theory of Poisson boundaries for
discrete quantum groups.

Fix a discrete quantum group $\widehat{\cG}$.
\begin{notation}
For every normal state $\phi\in \plon_*$, we define the
convolution operator
 $$P_{\phi}:\plon\recht\plon:P_{\phi}(a)=(\id \ot
 \phi)\widehat{\Delta}(a) \; .$$
\end{notation}

We are only interested in special states $\phi\in \plon$, motivated
by the following straightforward proposition. For every probability
measure $\mu$ on $\Irred(\cG)$, we set
\begin{equation*}
\psi_{\mu}=\sum_{x\in \Irred(\cG)}\mu(x)\psi_x \quad\text{and}\quad
P_\mu := P_{\psi_{\mu}} \; .
\end{equation*}
Recall that the states $\psi_x$ are defined in notation
\ref{not.states}. Note that we have a convolution product $\mu
* \nu$ on the measures on $\Irred(\cG)$, such that $\psi_{\mu * \nu}
= (\psi_\mu \ot \psi_\nu)\deh$.

\begin{proposition}
Let $\phi$ be a normal state on $\plon$. Then the following
conditions are equivalent. \begin{itemize}
\item $\phi$ has the form $\psi_{\mu}$ from some probability measure $\mu$ on
$\Irred(\cG)$.
\item The Markov operator $P_\phi$ preserves the center of $\ell^\infty(\cGh)$.
\item $\phi$ is invariant under the adjoint action of $\cG$ on $\ell^\infty(\cGh)$ \[\al_\cG :
  \ell^\infty(\cGh) \recht \ell^\infty(\cGh) \otvn L^\infty(\cG) : a \mapsto
  \mathbb{V}(a\ot 1)\mathbb{V}^*\;.\]
\end{itemize}
\end{proposition}

\begin{definition}[\cite{iz1}, Section 2.5]
Let $\mu$ be a probability measure on $\Irred(\cG)$. Set
 $$H^{\infty}(\widehat{\cG},\mu)=\{a \in \plon\mid\  P_{\mu}(a)=a\}
 \; .$$
Equipped with the product defined by
\begin{equation}\label{protput}
a \cdot b := \wlim{n \recht \infty} \frac{1}{n} \sum_{k=1}^n
P_\mu^k(ab) \; ,
\end{equation}
and the involution, norm and $\sigma$-weak topology inherited from
$\ell^\infty(\cGh)$,
 the space $H^{\infty}(\widehat{\cG},\mu)$ becomes
a von Neumann algebra that we call the \emph{Poisson boundary} of
$\cGh$ with respect to $\mu$.
\end{definition}

\begin{terminology}
A probability measure $\mu$ on $\Irred({\cG})$ is called
\emph{ge\-nerating} if there exists, for every $x\in \Irred({\cG})$,
an $n\geq 1$ such that $\mu^{* n}(x)\not=0$.
\end{terminology}

\begin{remark}
The restriction of the co-unit $\epsh$ yields a state on $\hon$,
called the \emph{harmonic state}. This state is faithful when $\mu$
is generating. In what follows, we always assume that $\mu$\ is
generating.
\end{remark}

\begin{definition} \label{def.actions}
Let $\mu$ be a generating measure on $\Irred(\cG)$. The Poisson
boundary $\hon$ comes equipped with two natural actions, one of
$\cG$ and one of $\cGh$:
\begin{align*}
\al_{\cG} & : \hon \recht \hon \otvn L^\infty(\cG) : \al_{\cG}(a) =
\mathbb{V}(a\ot 1)\mathbb{V}^* \; ,  \\
\al_{\cGh} & : \hon \recht \ell^\infty(\cGh) \otvn \hon :
\al_{\cGh}(a) =
  \deh(a) \; .
\end{align*}
\end{definition}
Note that $\al_{\cG}$ is the restriction of the adjoint action of
$\cG$ on $\ell^\infty(\cGh)$, while $\al_{\cGh}$ is nothing else
than the restriction of the comultiplication. The maps $\al_{\cG}$
and $\al_{\cGh}$ are well defined because of the following
equivariance formulae:
\begin{equation}\label{aactie}
(\id\ot P_{\mu})(\widehat{\Delta}(a))=\widehat{\Delta}(P_{\mu}(a))
\quad\text{and}\quad (P_{\mu}\ot
\id)(\al_\cG(a))=\al_\cG(P_{\mu}(a)) \; .
\end{equation}

\begin{remark}\label{kakaa} With the product defined by formula
$(\ref{protput})$, the mappings $\alpha_{\cG}$\ and $\alpha_{\cGh}$\ are multiplicative. This follows from the
equivariance formulae $(\ref{aactie})$. Hence $\alpha_\cG$\ and $\alpha_{\cGh}$\ are
actions on $\hon$. Because
$$(\epsh\ot \id)\alpha_{\cG}(a)=(\epsh\ot \id)(\mathbb{V}(a\ot
1)\mathbb{V}^*)=\epsh(a)1\;,$$ we see that $\epsh$\ is an invariant
state for the action $\alpha_{\cG}:\hon \recht \hon \ot
L^\infty(\cG)$.
\end{remark}

\begin{remark}
When $\widehat{\cG}$ is a discrete group, the action $\al_{\cG}$ is the trivial action
on $\plon$.
In general, the fixed point algebra of $\al_{\cG}$
is precisely the algebra of \emph{central harmonic elements}
$Z(\plon)\cap \hon$. Since the Markov operator $P_{\mu}$ preserves the
center $Z(\plon)$, the commutative von Neumann algebra $Z(\plon)\cap
\hon$ with state $\epsh$, is exactly the Poisson boundary for the random
walk on $\Irred(\cG)$ with transition probabilities
$p(x,y)$ and $n$-step transition probabilities $p_n(x,y)$
given by
\begin{equation}\label{eq.nstep}
p_x p(x,y) = p_x P_\mu(p_y) \; , \quad p_x p_n(x,y) = p_x P_\mu^n(p_y)
\; .
\end{equation}
Note that $p_n(e,y) = \mu^{* n}(y) = \psi_\mu^{*
n}(p_y)$.

So, the action $\al_{\cG}$ is \emph{ergodic} if and
only if there are no non-trivial central harmonic elements.
\end{remark}

\section{The Martin boundary of a discrete quantum group}

The Martin boundary and the Martin compactification of a discrete
quantum group have been defined by Neshveyev and Tuset in
\cite{NT}.
Fix a discrete quantum group $\cGh$ and
a probability measure $\mu$ on $\Irred(\cG)$. We have an
associated Markov operator $P_\mu$ and a classical random walk on $\Irred(\cG)$ with $n$-step transition probabilities
given by \eqref{eq.nstep}.

\begin{definition}
The probability measure $\mu$ on $\Irred(\cG)$ is said to be \emph{transient}
if $\sum_{n=0}^\infty p_n(x,y) < \infty$ for all $x,y \in \Irred(\cG)$.
\end{definition}

We suppose throughout that $\mu$ is a
generating measure and that $\mu$ is transient.

Denote by $c_c(\cGh) \subset c_0(\cGh)$ the algebraic direct sum of
the algebras $\cL(H_x)$. We define, for $a \in c_c(\cGh)$,
$$G_\mu(a) = \sum_{n=0}^\infty P_\mu^n(a) \; .$$
Observe that usually $G_\mu(a)$ is unbounded, but it makes sense in
the multiplier algebra of $c_c(\cGh)$, i.e.\ $G_\mu(a)p_x \in \cL(H_x)$
makes sense for every $x \in \Irred(\cG)$ because $\mu$ is transient. Moreover, $G_\mu(p_\veps)$ is strictly positive and
central. This allows to define the Martin kernel as follows.

Whenever $\mu$ is a measure on $\Irred(\cG)$, we use the notation
$\overline{\mu}$ to denote the measure given by $\overline{\mu}(x) = \mu(\overline{x})$.

\begin{definition}[Defs.\ 3.1 and 3.2 in \cite{NT}]
Define
$$K_\mu : c_c(\cGh) \recht \ell^\infty(\cGh) : K_\mu(a) = G_\mu(a)
G_\mu(p_\veps)^{-1} \; .$$
Define the \emph{Martin compactification} $\widetilde{A}_{\mu}$ as the
\cst-subalgebra of $\ell^\infty(\cGh)$ generated by
$K_{\overline{\mu}}(c_c(\cGh))$ and $c_0(\cGh)$. Define the \emph{Martin
  boundary} $A_{\mu}$ as the quotient $\widetilde{A}_{\mu} / c_0(\cGh)$.
\end{definition}

By Theorem 3.5 in \cite{NT}, the adjoint action $\al_\cG$ and the
comultiplication $\deh$ define, by restriction
\begin{equation}\label{eq.actionsMartin1}
\alpha_\cG : \widetilde{A}_\mu \recht \widetilde{A}_\mu \ot C(\cG) \quad\text{and}\quad
\alpha_{\cGh} : \widetilde{A}_\mu \recht \M(c_0(\cGh) \ot \widetilde{A}_\mu) \; .
\end{equation}
By passing to the quotient, we get the following actions on the Martin boundary.
\begin{equation}\label{eq.actionsMartin}
\gamma_\cG : A_\mu \recht A_\mu \ot C(\cG) \quad\text{and}\quad
\gamma_{\cGh} : A_\mu \recht \M(c_0(\cGh) \ot A_\mu) \; .
\end{equation}

\begin{remark} The actions $\alpha_\cG$\ and $\gamma_\cG $\ are {\it reduced}.
\end{remark}

\section{The correspondence between the actions of monoidally equivalent quantum groups}

In this section, we prove that there is a {\it bijective correspondence} between actions of monoidally
equivalent compact quantum groups. Moreover, this correspondence preserves the spectral properties of the actions.

\subsection{Construction of the bijective correspondence}

\begin{notation}
Consider two compact quantum groups $\cG_1$ and $\cG_2$ and a *-algebra $\cB$ on which there exist two commuting actions $\sde_1: \cB\to\cC(\cG_1)\otalg \cB$ and $\sde_2:\cB\to\cB\otalg\cC(\cG_2)$. Given an action $\alpha: \cD\to\cD\otalg\cC(\cG_1)$, we define
$$\cD\boxalg^{\alpha}\cB:=\{a\in D\otalg\cB\mid (\alpha\ot\id)(a)=(\id\ot\sde_1)(a)\}\;.$$
When we consider everything on the von Neumann algebraic level, we denote in the same way
$$D\boxtimes^{\alpha}B:=\{a\in D\otvn B\mid (\alpha\ot\id)(a)=(\id\ot\sde_1)(a)\}\;.$$
Also, if $D$\ is a $C^*$-algebra with reduced action $\alpha$\ on it, then we denote
$D\boxtimes_{red}^{\alpha}B$\ as the norm closure of $\cD\boxalg^{\alpha}\cB$\ in $D\boxtimes^{\alpha}B$.
\end{notation}

\begin{lemma}
The restriction of $\id\ot\sde_2$ to $\cD\boxalg^{\alpha}\cB$ gives an action of $\cG_2$ on $\cD\boxalg^{\alpha}\cB$. We denote this action by $\id\boxtimes \sde_2$.
\end{lemma}

\begin{proof}
From the following easy calculation, one can see that $\cD\boxalg^{\alpha_1}\cB$
is invariant under the action $\id\ot\sde_2$.
\begin{align*}
(\alpha\ot\id\ot\id)(\id\ot\sde_2)(a)&=(\id\ot\id\ot\sde_2)(\alpha\ot
\id)(a)=(\id\ot\id\ot\sde_2)(\id\ot\sde_1)(a)\\&=(\id\ot\sde_1\ot\id)(\id\ot\sde_2)(a)
\end{align*}
The last step is valid because $\sde_1$ and $\sde_2$ commute. Hence $\id\boxtimes \sde_2$ is a well defined action of
$\cG_2$ on $\cD\boxalg^{\alpha}\cB$.
\end{proof}

Consider two monoidally equivalent compact quantum groups $\cG_1$
and $\cG_2$ and a \cst-algebra $D_1$. Suppose we have an action
$\alpha_1:D_1\recht D_1\ot C(\cG_1)$. As we stated in remark
\ref{hopf} and remark \ref{prop.densesub} , the underlying
Hopf*-algebra action carries all the relevant information. This
means that we can work with this underlying Hopf $^*$-algebra
action $\alpha_1:\mathcal{D}_1\recht \mathcal{D}_1\otalg
\mathcal{C}(\cG_1)$. Consider a monoidal equivalence
$\varphi:\cG_2\recht \cG_1$. Note that we have exchanged the roles
of $\cG_1$ and $\cG_2$. This will turn out to be more convenient
in what follows. From theorem \ref{construction},  we get a link
algebra $\cB$, unitaries $X^x\in\cL(H_x,H_{\vphi(x)})\otalg \cB$
and two commuting ergodic actions
\[\sde_1: \cB\to \mathcal{C}(\cG_1)\otalg \cB \qquad\text{   and   }\qquad\sde_2:\cB\to
\cB\otalg \mathcal{C}(\cG_2)\;.\] given by
\begin{equation}\label{qqsszz}
(\id\ot \sde_1)(X^x)=U^{\varphi(x)}_{12}X^x_{13}\qquad\text{   and
}\qquad (\id\ot \sde_2)(X^x)=X^x_{12}U^x_{13}
\end{equation}

The following theorem enables us to construct an action of $\cG_2$
with the same spectral structure as $\alpha_1$.

\begin{theorem}\label{hierse}
The action $\alpha_2:=\id\boxtimes \sde_2$\ of $\cG_2$\ on $\cD_2:=\cD_1\boxalg^{\alpha_1}\cB$\
has the following properties:
\begin{itemize}
\item $a \mapsto a \ot 1$ is a $^*$-isomorphism between the fixed point algebras of $\alpha_1$\ and $\alpha_2$.

\item The map $T_x:K_{\vphi(x)}\to K_x: v\mapsto v_{12}X^x_{13}$ is a bimodular isomorphism between the spectral
subspaces of $\alpha_1$ and $\alpha_2$. Moreover, $T$\ is a
unitary element of $\cL(K_{\vphi(x)},K_x)$ for the inner products
$\langle\cdot,\cdot\rangle_l$ and  $\langle\cdot,\cdot\rangle_r$
defined by formulae (\ref{hil1}) and (\ref{hil2}).

\item The set $(T_x)_{x\in\Irred(\cG_2)}$ respects the monoidal
structure in the sense that for $x,y,z\in\Irred(\cG_2)$ and
$V\in\Mor(x\ot y,z)$
\[T_x(X)_{13}T_y(Y)_{23}(V\ot 1)=T_z(X_{13}Y_{23}(\vphi(V)\ot 1))\]
for all $X\in K_{\vphi(x)}$, $Y\in K_{\vphi(y)}$.

\item Suppose that $\alpha_1$\ is an ergodic action. Then the action $\alpha_2$ as defined above is also
 an ergodic action. Moreover
 for all $x\in \operatorname{Irred}(\cG_2)$, $\mult_q(x)=\mult_q(\vphi(x))$.

\end{itemize}
\end{theorem}

\begin{proof}

Suppose that $\alpha_2(a)= a\ot 1$ for $a\in \cD_2$. This means
that $(\id\ot\sde_2)(x)=x\ot 1$. By ergodicity of $\sde_2$ there
exists a $b\in\cD_1$ such that $a=b\ot 1$. But because
$(\alpha_1\ot\id)(a)=(\id\ot\sde_1)(a)$, it follows that
$b\in\cD_1^{\alpha_1}$. So the map
$\cD_1^{\alpha_1}\to\cD_2^{\alpha_2}:b\mapsto b\ot 1$\ is a
*-isomorphism.

We now prove that the spectral subspaces of $\alpha_1$ and
$\alpha_2$ are isomorphic as $\cD_1^{\alpha_1}$-bimodules. Denote
by $K_{\vphi(x)}$ and $K_x$ the spectral subspaces of respectively
$\alpha_1$ and $\alpha_2$ for the representation $\vphi(x)$,
respectively $x$. From remark \ref{Hilbert}, we know that the
spectral subspaces have a natural bimodule structure over the
fixed point algebra. We claim that the map \[T:K_{\vphi(x)}\to
K_x: v\mapsto v_{12}X^x_{13}\] is the bimodule isomorphism we are
looking for. If $v\in K_{\vphi(x)}$, then
\[
(\id\ot\alpha_1\ot \id)T(v)=v_{12}U^{\vphi(x)}_{13}X^x_{14}
=(\id\ot\id\ot\sde_1)T(v)
\]
by definition \ref{def.spectral} of the spectral subspace
$K_{\vphi(x)}$
 and the properties of $X^x$, so
$T(v)\in\overline{H_x}\ot \cD_2$. Moreover, it is obvious that
\[(\id\ot\alpha_2)T(v)=(\id\ot\id\ot\sde_2)T(v)=v_{12}X^x_{13}U^x_{14}\;,\]
which means $T(v)\in K_x$. The $\cD_1^{\alpha_1}$-bilinearity of $T$
is clear. Consider now the spectral subspaces $K_x$\ and
$K_{\varphi(x)}$\ as equipped with the left inner product as in
(\ref{hil1}). We show that $T$\ is a unitary element of
$\cL(K_{\vphi(x)},K_x)$ for this inner product and obtain in this
way that $T$\ actually gives an isomorphism between $K_{\vphi(x)}$
and $K_x$. Consider the map $S:K_x\to \overline{H_{\vphi(x)}}\ot
\cD_1\ot \cB:w\mapsto w(X^x_{13})^*$. If $w\in K_{x}$, then
\begin{align*}
(\id\ot\id\ot
\sde_2)S(w)&=(\id\ot\alpha_2)(w)(\id\ot\id\ot\sde_2)(X^x_{13})^*=(w\ot
1)U^x_{14}(U^x_{14})^*(X^x_{13})^*=S(w)\ot 1\;,
\end{align*}
So, by ergodicity of $\sde_2$, we may conclude that $S(w)\in
\overline{H_{\vphi(x)}}\ot\cD_1\ot \mathbb{C}$.

Because $w$ has its second leg in $\cD_2$, we get that
\begin{align*}
(\id\ot\alpha_1\ot\id)S(w)&=(\id\ot\id\ot\sde_1)(w)(X^x_{14})^*=(\id\ot\id\ot\sde_1)(w(X^x_{13})^*)U^{\varphi(x)}_{13}\\
&=(\id\ot\id\ot\sde_1)(S(w))U^{\varphi(x)}_{13}\;,
\end{align*}

 But we just proved that the third leg of $S(w)$ is scalar, so the last expression is nothing else
than $(S(w)\ot 1)U^{\vphi(x)}_{13}$. Thus, by the definition of
$K_{\varphi(x)}$, we get that $S: K_x\to K_{\vphi(x)}\ot \C$.

For every  $v\in K_{\vphi(x)}$ and $w\in K_x$, we have that
\[
\langle T(v),w\rangle_l=T(v)w^*=v_{12}X^x_{13}w^*=v_{12}S(w)^*
=\langle v, S(w)\rangle_l\;.\]  So,  $S$ is actually the adjoint $T^*$\ of $T$ in the
sense of Hilbert \cst-modules. Moreover, it is trivial that
$T^*T=1=TT^*$. Hence $T\in\cL(K_{\vphi(x)},K_x)$ is unitary for
$\langle\cdot,\cdot\rangle_l$.

Next, we show that $T$ is also a unitary element of
$\mathcal{L}(K_{\varphi(x)},K_x)$\ for the right
Hilbert C$^*$-module structure given by (\ref{hil2}). From
proposition $3.5$\ of \cite{lance}, it suffices to show that $T$ is
isometric and surjective. The surjectivity follows from above. We
use the orthogonality relations (\ref{eq.orthogonality}) and
(\ref{eq.orthogonalitymon}) for $U^x$ and $X^x$\ to prove that $T$\
is indeed an isometry.

First notice that the conditional expectation
$E_2:\cD_2\to\cD_2^{\alpha_2}$ is nothing else than the map
$a\mapsto (\id\ot\om)(a)\ot 1$, where $\om $\ is the invariant state
for $\sde_1$ and $\sde_2$. Indeed, for $a\in \cD_2$,
\[E_2(a)=(\id\ot h_2)\alpha_2(a)=(\id\ot\id\ot h_2)(\id\ot\sde_2)(a)=(\id\ot\om)(a)\ot
1\;.\]

Consider now $v\in K_{\vphi(x)}$. On the one hand, we have that
\begin{align*}
1\ot \langle T(v), T(v)\rangle_{r}&=(\id\ot
E_2)((X^x_{13})^*v^*_{12}v_{12}X^x_{13})=(\id\ot\id\ot\om)
((X^x_{13})^*v^*_{12}v_{12}X^x_{13})\ot 1\\
&=\frac{1}{\dimq(x)}(1\ot(\Tr \ot\id)((Q^{-1}_{\vphi(x)}\ot
1)v^*v)\ot 1)
\end{align*}

because of the orthogonality relations for $X^x$.

On the other hand
\begin{align*}
1\ot \langle v,v \rangle_{r}&=(\id\ot E_1)(v^*v)=(\id\ot\id\ot
h_1)(\id\ot\alpha_1)(v^*v)=(\id\ot\id\ot
h_1)((U^x_{13})^*v^*_{12}v_{12}U^x_{13})\\&=\frac{1}{\dimq(x)}(1\ot(\Tr
\ot\id)((Q^{-1}_{\vphi(x)}\ot 1)v^*v))\;,
\end{align*}
where in the last step we used the orthogonality relations for
$U^x$. Considering the map $\mathcal{D}_1^{\alpha_1}\recht
\mathcal{D}_2^{\alpha_2}:a\mapsto a\ot 1$, the calculations above
show that $T$ is indeed isometric and hence unitary.

We now show that $(T_x)_{x\in\Irred(\cG_2)}$ preserves the monoidal
structure. Take $v\in K_{\vphi(x)}$, $w\in K_{\vphi(y)}$ and
$V\in\Mor(x\ot y,z)$. We calculate that
\begin{align*}
T_x(v)_{13}T_y(w)_{23}(V\ot 1)&=v_{13}X^x_{14}w_{23}X^y_{24}(V\ot
1)\\&=v_{13}w_{23}(\vphi(V)\ot
1)X^z_{13}=T_z(v_{13}w_{23}(\vphi(V)\ot 1))\;,
\end{align*}
which proves the statement.

Finally, we prove the fourth part of the theorem.  Recall
 the operators from formula (\ref{teees}).
$$R_{\varphi(x)}:v\mapsto (\varphi(t)^* \ot 1)(1 \ot v^*)\qquad \mbox{and}\qquad L_{\varphi(x)}=R_{\varphi(x)}^*R_{\varphi(x)}$$
 with $v\in K_{\varphi(x)}$
and
$$R_{x}:w\mapsto (t^* \ot 1)(1 \ot w^*)\qquad \mbox{and}\qquad L_{x}=R_{x}^*R_{x}$$
with $w\in K_{x}$. Then
\begin{align*}
\langle v, L_{\varphi(x)}w\rangle 1&=\langle
R_{\varphi(x)}v,R_{\varphi(x)}w\rangle 1= \langle (\varphi(t)^*\ot
1)(1\ot v^*),(\varphi(t)^*\ot 1)(1\ot w^*)\rangle\\
&= (\varphi(t)^*\ot 1)(1\ot w^*v)(\varphi(t)\ot 1)
\end{align*}
where $v,w\in K_{\varphi(x)}$.
 Remember the isomorphism $T_x:K_{\vphi(x)}\to
K_x: v\mapsto v_{12}X^x_{13}$. Then:
\begin{eqnarray*}
\langle v_{12}X^x_{13}, L_x w_{12}X^x_{13}\rangle 1&=&\langle
R_{x} v_{12}X^x_{13},R_{x} w_{12}X^x_{13}\rangle 1
= (t^*\ot 1)(1\ot (w_{12}X^x_{13})^*(v_{12}X^x_{13}))(t\ot 1)\\
&=& (t^*\ot 1\ot 1) ((X^x_{24})^*w_{23}^*v_{23}X^x_{24})(t\ot 1\ot
1)\\ &=& (\varphi(t)^*\ot 1\ot 1)
(X^{\overline{x}}_{14}w_{23}^*v_{23}(X^{\overline{x}}_{14})^*)(\varphi(t)\ot
1\ot 1)\\
&=& (\varphi(t)^*\ot 1\ot 1)(w_{23}^*v_{23})(\varphi(t)\ot 1\ot 1)
= (\varphi(t)^*\ot 1)(1\ot w^*v)(\varphi(t)\ot 1)\ot 1\;.
\end{eqnarray*}
In this calculation, we have used that
$X^{\overline{x}}_{13}X^x_{23}(t\ot 1)=\varphi(t)\ot 1$. This
follows from the fact that $t\in \Mor(\overline{x}\ot
x,\varepsilon)$. Again considering the map
$\mathcal{D}_1^{\alpha_1}\recht \mathcal{D}_2^{\alpha_2}:a\mapsto
a\ot 1$, we get that $T_x$ intertwines $L_x$\ and $L_{\varphi(x)}$.
It follows trivially from the definition \ref{quam} of quantum
multiplicity that both quantum multiplicities are the same. This
completes the proof of the theorem.
\end{proof}

\begin{remark}\label{vonnie}
It seems that the statement of the theorem cannot
immediately be formulated on the \cst-algebraic level. If we
define $D_2=\{a\in D_1\ot B\mid
(\alpha_1\ot\id)(a)=\id\ot\sde_1)(a)\}$, as we did before for the algebraic and the von Neumann algebraic case, it is not clear that
$\alpha_2=\id\ot\sde_2$ goes into $D_2\ot C(\cG_2)$.

However, for von Neumann algebras, there is no problem. Suppose
that $\alpha_1:D_1\to D_1\otvn L^{\infty}(\cG)$ is a von Neumann
algebraic action and take the notations as before, where we now
take the von Neumann algebraic link-algebra $B=(\cB,\om)''$. Here
it does hold that $(\id\ot \sde_2)(D_1\boxtimes^{\alpha_1}B)\subset D_1\boxtimes^{\alpha_1}B\otvn
L^{\infty}(\cG_2)$, as we only need to check that
$(\id\ot\id\ot\mu)(\id\ot \sde_2)(a)\in D_1\boxtimes^{\alpha_1}B$ for all $a\in D_1\boxtimes^{\alpha_1}B$ and $\mu\in
(L^{\infty}(\cG_2))_{*}$. For \cst-algebras, this argument is not valid.
\end{remark}
\label{von neumann level}

\textbf{Claim:} The algebra $\cD_2:=\cD_1\boxalg^{\alpha_1}\cB$\ as defined in theorem
\ref{hierse} is precisely the spectral subalgebra of
$(D_1\boxtimes^{\alpha_1}B,\id\boxtimes \sde_2)$.

\begin{proof}
Denote by $\widetilde{\cD_2}$ the spectral subalgebra of
$(D_1\boxtimes^{\alpha_1}B,\id\boxtimes \sde_2)$. It is clear that $\cD_1\boxalg^{\alpha_1}\cB\subset \widetilde{\cD_2}$.

On the other hand, \[\widetilde{\cD_2}=\langle(\id\ot
h)((\id\boxtimes \sde_2)(a)(1\ot b))\mid a\in D_1\boxtimes^{\alpha_1}B, b\in \cC(\cG_2)_x, x\in
\Irred(\cG)\rangle\;.\] Because the elements of $\widetilde{\cD_2}$
of course sit in $D_1\boxtimes^{\alpha_1}B$, it is sufficient to prove that
$\widetilde{\cD_2}\subset \cD_1\otalg\cB$.

If $a\in D_1\boxtimes^{\alpha_1}B$, $b\in \cC(\cG_2)_x$, $x\in \Irred(\cG_2)$, $(\id\ot
h)((\id\boxtimes \sde_2)(a)(1\ot b))$ belongs to the strongly closed linear span
of \[\{(\id\ot\id\ot h)((a\ot \sde_2(d))(1\ot 1\ot b))\mid a\in
D_1,d\in\cB,b\in \cC(\cG_2)_x\}\subset D_1\otalg \cB_x\;.\] Since
$\cB_x$ is finite dimensional, $D_1\otalg \cB_x$ is already strongly
closed in $D_1\otvn B$. Hence $c:=(\id\ot h)((\id\boxtimes \sde_2)(a)(1\ot b))\in
D_1\otalg \cB_x$ for all $a\in D_1\boxtimes^{\alpha_1}\cB$ and $b\in \cC(\cG_2)_x$. On the
other hand, $(\alpha_1\ot\id)(c)=(\id\ot\sde_1)(c)$, which implies
that $c\in (\cD_1)_x\otalg \cB_x\subset \cD_1\otalg\cB$. This ends
the proof.
\end{proof}

We can start from the comultiplication on $\cG_1$ and apply the
above construction. It is not surprising that we end up with the
link algebra and the action $\sde_2$.

\begin{proposition}\label{link algebra}
Consider two monoidally equivalent compact quantum groups $\cG_1$\
and $\cG_2$. Then there exists a $^*$-isomorphism between $\cC(\cG_1)\boxalg^{\de_1}\cB$\ and
the link algebra $\cB$.
Moreover, this $^*$-isomorphism intertwines the action $\id\boxtimes \sde_2$\
with the action $\sde_2$.
\end{proposition}
\begin{proof}
 We claim that $\sde_1:\cB\to \cC(\cG_1)\boxalg^{\de_1}\cB$ is the desired $^*$-isomorphism. From the definition of
$\sde_1$, it follows that $\sde_1: \cB\to \mathcal{C}(\cG_1)\otalg
\cB$\
 is an injective $^*$-homomorphism. The image of $\sde_1$\ is contained in $\cC(\cG_1)\boxalg^{\de_1}\cB$\ because $\sde_1$\ is an action.
   Moreover, if $a\in \cC(\cG_1)\boxalg^{\de_1}\cB$ and $\eps_1$ is the co-unit on $\mathcal{C}(\cG_1)$, then
\[\sde_1((\eps_1\ot\id)(a))=(\eps_1\ot\id\ot\id)(\de_1\ot\id)(a)=a\]
which means that $\sde_1(\cB)=\cC(\cG_1)\boxalg^{\de_1}\cB$. So $\sde_1$\ is also
surjective.

Because $\sde_1$ and $\sde_2$ commute, it is clear that this
$^*$-isomorphism intertwines the actions $\sde_2$\ and $\id\boxtimes \sde_2$.
\end{proof}

Now we consider the inverse monoidal equivalence
$\vphi^{-1}:\cG_1\to \cG_2$. According to theorem
\ref{construction}, we obtain the link algebra $\tilde{\cB}$\
generated by the coefficients of unitary elements $Y^x\in
\cL(H_{\vphi(x)},H_x)\otalg \tilde{\cB}$\ and two commu\-ting
ergodic actions
\[\gamma_2:
\tilde{\cB}\to \mathcal{C}(\cG_2)\otalg\tilde{\cB}\qquad \text{and}\qquad \gamma_1:\tilde{\cB}\to  \tilde{\cB}\otalg\mathcal{C}(\cG_1)\] with
\begin{equation}\label{vvxxww}
(\id\ot \gamma_2)(Y^x)=U^x_{12}Y^x_{13} \qquad
\text{and}\qquad (\id\ot \gamma_1)(Y^x)=Y^x_{12}U^{\varphi{(x)}}_{13}\;.
\end{equation}
 Denote by
$\tilde{\om}$ the invariant state on $\tilde{\cB}$. Then we get the
following proposition.

\begin{proposition}\label{quantumgroup}
Consider two monoidally equivalent compact quantum groups $\cG_1$\
and $\cG_2$. Then we obtain a $*$-isomorphism
$\pi:\mathcal{C}(\cG_2)\to \tilde{\cB}\boxalg^{\gamma_1}\cB$. This $^*$-isomorphism intertwines
the comultiplication  $\Delta_2$\ with the action $\id\boxtimes \sde_2$.
\end{proposition}

\begin{proof}
 We define the
linear map $\pi:\mathcal{C}(\cG_2)\to \tilde{\cB}\ot \cB$ where
$(\id\ot\pi)(U^x)=Y^{x}_{12}X^x_{13}$. Because
\begin{align*}
(\id\ot \gamma_1\ot
\id)(Y^x_{12}X^x_{13})=Y^x_{12}U^{\varphi(x)}_{13}X^x_{14}=(\id\ot\id
\ot\sde_1)(Y^x_{12}X^x_{13})\;,
\end{align*}
the image of $\pi$ lies in $\tilde{\cB}\boxalg^{\gamma_1}\cB$.

Consider $x,y,z \in\Irred(\cG_2) $ and take $T\in\Mor(x\ot y, z)$.
The multiplicativity of $\pi$ follows from the following
calculation:
\begin{align*}
(\id\ot\pi)(U^x_{13}U^y_{23}(T\ot1))&=(\id\ot\pi)((T\ot 1)U^z)
=(T\ot 1\ot 1)(Y^z_{12}X^z_{13})\\&=Y^x_{13}Y^y_{23}(\vphi(T)\ot
1\ot 1)X^z_{13}=Y^x_{13}Y^y_{23}X^x_{14}X^y_{24}(T\ot 1\ot
1)\\&=(\id\ot\pi)(U^x)_{134}(\id\ot\pi)(U^y)_{234}(T\ot 1\ot 1)\;.
\end{align*}

Take now $t_x\in\Mor (x\ot\ox,\veps)$. Because
\[U^x_{13}U^{\ox}_{23}(t_x\ot 1)=t_x\ot 1\;,\] it follows that

\begin{align*}
(\id\ot\id\ot\pi)({U^{x}_{13}}^*(t_x\ot
1))&=(\id\ot\id\ot\pi)(U^{\ox}_{23}(t_x\ot
1))=Y^{\ox}_{23}X^{\ox}_{24}(t_x\ot 1\ot
1)\\&=Y^{\ox}_{23}{X^x_{14}}^*(\vphi(t_x)\ot 1\ot
1)={X^x_{14}}^*{Y^x_{13}}^*(t_x\ot 1\ot 1)\;.
\end{align*}
This proves that $\pi$ also passes trough the involution, so it is a
*-homo\-morphism. We now show that this map is the desired
*-isomorphism.

First we prove the injectivity. It is easy to show that
$(\tilde{\om}\ot\om)\pi=h$, with $h$ the Haar measure of
$\mathcal{C}(\cG_2)$. Suppose now that for an $a\in
\mathcal{C}(\cG_2)$, $\pi(a)=0$. Then also $\pi(a^*a)=0$, which
means that also $h(a^*a)=0$. But $h$ is faithful on
$\mathcal{C}(\cG_2)$, so $a=0$.

To prove the surjectivity of $\pi$, we have to take a closer look at
the elements of $\tilde{\cB}\boxalg^{\gamma_1}\cB$. From definition \ref{def.densesub},
we get that
$$\cB=\bigoplus_{x\in\mbox{Irred}(\cG_1)}\cB_x\qquad
\mbox{and}\qquad
\tilde{\cB}=\bigoplus_{x\in\mbox{Irred}(\cG_1)}{\tilde{\cB}}_x$$
where $\sde_1(\cB_x)\subseteq \mathcal{C}(\cG_1)_x\otalg \cB_x$ and
$\gamma_1({\tilde{\cB}}_x)\subseteq {\tilde{\cB}}_x\otalg
\mathcal{C}(\cG_1)_x$. Suppose $b\in \tilde{\cB}\boxalg^{\gamma_1}\cB$. We claim that
$b\in\bigoplus_{x\in \mbox{Irred}(\cG_1)}\tilde{\cB}_x\otalg\cB_x$.
First notice that $$\tilde{\cB}\otalg\cB=\bigoplus_{x,y\in
\mbox{Irred}(\cG_1)}\tilde{\cB}_x\otalg \cB_y\;.$$ So $b=\sum
b_{xy}$ with $b_{xy}\in \tilde{\cB}_x\otalg \cB_y$. Since
$\sde_1(\cB_y)\subseteq \mathcal{C}(\cG_1)_x\otalg \cB_y$\ and
$\gamma_1({\tilde{\cB}}_x)\subseteq {\tilde{\cB}}_x\otalg
\mathcal{C}(\cG_1)_x$, it follows that $b_{xy}=0$ if $x\neq y$.

So we only need to prove that
$\pi(\mathcal{C}(\cG_2)_x)=\tilde{\cB}_x\boxalg^{\gamma_1}\cB_x$,

Therefore, remember the formulas $(\id\ot
\gamma_1)(Y^x)=Y^x_{12}U^{\varphi{(x)}}_{13}$\ and $(\id\ot
\sde_1)(X^x)=U^{\varphi{(x)}}_{12} X^x_{13}$\ where $Y^x\in
\cL(H_{\varphi(x)},H_x)\otalg \tilde{\cB}_x$\ and
$X^x\in\cL(H_x,H_{\varphi(x)})\otalg \cB_x$.

We know that a basis of $\cB_x$\ (resp. $\tilde{\cB}_x$) is given
respectively by elements of the form $(\om_{\widetilde{e_
{k_x}},e_{l_x}}\ot\id)(X^x)$\ and $(\om_{e_{l_x},\widetilde{e_
{k_x}}}\ot\id)(Y^x)$\ with $\widetilde{e_ {k_x}}, k_x\in \{1,\ldots,
\dim(\varphi(x)) \}$\  an orthonormal basis in $H_{\varphi(x)}$\ and
${e_ {l_x}}, l_x\in$\linebreak$\{1,\ldots,\dim(x) \}$\ an
orthonormal basis in $H_x$. Denote $(\om_{\widetilde{e_
{k_x}},e_{l_x}}\ot\id)(X^x):=b_{k_x,l_x}$\ and
$(\om_{e_{l_x},\widetilde{e_
{k_x}}}\ot\id)(Y^x):=\tilde{b}_{l_x,k_x}$. We also have a basis for
$\mathcal{C}(\cG_2)_x$\ given by $(\om_{{e_
{k_x}},e_{l_x}}\ot\id)(U^x)$, again with ${e_ {k_x}}, k_x\in
\{1,\ldots, \dim(x) \}$\ an orthonormal basis in $H_x$. Denote by
$(\om_{{e_ {k_x}},e_{l_x}}\ot\id)(U^x):=u_{k_x,l_x}$. In the
following, we drop the subscript $x$. With these notations, we get
that
$$\gamma_1(\tb_{kl})=\sum_{p=1}^{\dim (\varphi(x))}\tb_{kp}\ot
u_{pl}\qquad \mbox{and}\qquad \sde_1(b_{ij})=\sum_{q=1}^{\dim
(\varphi(x))}u_{iq}\ot x_{qj}$$ and
$$\pi(u_{kj})=\sum_{l=1}^{\dim (\varphi(x))}\tb_{kl}\ot b_{lj}\;.$$
An arbitrary element of $\cB_x\otalg \tilde{\cB}_x$\ has the form
$$a=\sum_{klij}\lambda^{kj}_{li} \tb_{kl}\ot b_{ij}\;.$$ We get that
$$(\gamma_1\ot \id)(a)=\sum_{klijp}\lambda^{kj}_{li} \tb_{kp}\ot u_{pl}\ot b_{ij}$$
equals
$$(\id \ot \sde_1)(a)=\sum_{klijq}\lambda^{kj}_{li} \tb_{kl}\ot u_{iq}\ot b_{qj}\;.$$
From this equality, we immediately see that $\lambda^{kj}_{li}=0$\
if $l\neq i$. We have also that $\lambda^{kj}_{ll} =
\lambda^{kj}_{pp}$\ for all $k,l,j,p$.
 Indeed, the above equality provides us with the following equalities:
$$\sum_{lq}\lambda^{kj}_{ll}\tb_{kl}\ot u_{lq}\ot b_{qj}=\sum_{lp}\lambda^{kj}_{ll}\tb_{kp}\ot u_{pl}\ot b_{lj}$$
for every $k,j$. This can only happen when
$\lambda^{kj}_{ll}=\lambda^{kj}$\ for every $l\in \{1, \ldots
\dim(H_{\varphi(x)}) \}$. So $a\in \tilde{\cB}_x\boxalg^{\gamma_1}\cB_x$\ has the form
$$a=\sum_{kj}\lambda^{kj}\left(\sum_{l=1}^{\dim(\varphi(x))}\tb_{kl}\ot b_{lj}\right)$$
which is a linear combination of the $\pi(u_{kj})$. This proves the
surjectivity of $\pi$.

Moreover
\begin{align*}
(\id\ot\pi\ot\id)(\id\ot\de_2)(U^x)&=(\id\ot\pi\ot\id)(U^x_{12}U^x_{13})
=Y^x_{12}X^x_{13}U^x_{14}=(\id\ot(\id\ot \sde_2)\pi)(U^x)\;,
\end{align*}
so the action $\id\boxtimes \sde_2$ indeed corresponds to the comultiplication
on $\cG_2$.
\end{proof}

A combination of the two previous propositions now enables us to
prove the reversibility of our construction.

\begin{proposition}\label{omgekeerd}
Consider two monoidally equivalent compact quantum groups $\cG_1$\
and $\cG_2$\ and suppose also that $\alpha_1: \cD_1\recht \cD_1\ot \cC(\cG)$\ is an action. Then $\cD_1$\ and
$(\cD_1\boxalg^{\alpha_1}\cB)\boxalg^{\alpha_2} \tilde{\cB}$\ are $*$-isomorphic. Moreover,  this $^*$-isomorphism intertwines the actions $\alpha_1$\ and $\alpha_1'=(\id\ot\id\ot \gamma_1)\mid_{(\cD_1\boxalg^{\alpha_1}\cB)\boxalg^{\alpha_2} \tilde{\cB}}$.

\end{proposition}
\begin{proof}
 In exactly the same way as in proposition \ref{quantumgroup}, we can prove that $\mathcal{C}(\cG_1)$\ is $^*$-isomorphic to
$\cB\boxalg^{\delta_2}\tilde{\cB}$. In this case, the
$^*$-isomorphism is given by $\pi:\mathcal{C}(\cG_1)\to \cB\boxalg^{\delta_2}\tilde{\cB}$ where
$(\id\ot\pi)(U^{\varphi(x)})=X^{x}_{12}Y^x_{13}$.  Also in the same
way, we can prove that $\pi$\ intertwines the actions $\sde_1\ot
\id\mid_{\cB\boxalg^{\delta_2}\tilde{\cB}}$\ and $\Delta_1$. From this, we get that $(\cD_1\boxalg^{\alpha_1}\cB)\boxalg^{\alpha_2} \tilde{\cB}$ is
isomorphic to \[\cD_1'':=\{a\in\cD_1\otalg \mathcal{C}(\cG_1)\mid
(\alpha_1\ot\id)(a)=(\id\ot\de_1)(a)\}\;.\] From the calculation
\begin{align*}
(\id\ot\pi\ot\id)(\id\ot\de_1)(U^{\varphi(x)})&=(\id\ot\pi\ot\id)(U^{\varphi(x)}_{12}U^{\varphi(x)}_{13})=X^x_{12}Y^x_{13}U^{\varphi(x)}_{14}
\\&=(\id\ot\id\ot \gamma_1)(\id\ot\pi)(U^{\varphi(x)})\;,
\end{align*}
it follows that $(\id\ot\id\ot \gamma_1)\mid_{(\cD_1\boxalg^{\alpha_1}\cB)\boxalg^{\alpha_2} \tilde{\cB}}$\ is
equivalent with $\id\ot \de_1\mid_{\cD_1''}$\ under the
$^*$-isomorphism $\id\ot\pi$.

In the same way as in proposition \ref{link algebra}, we see that
$\alpha_1:\cD_1\to\cD_1''$ is a $^*$-isomorphism. It is obvious that
this $^*$-isomorphism  intertwines the actions $\id\ot
\de_1\mid_{\cD_1''}$\ and $\alpha_1$. Thus we get that
\begin{equation}\label{nogeffee}
(\id\ot \pi)\circ \alpha_1:\cD_1\recht (\cD_1\boxalg^{\alpha_1}\cB)\boxalg^{\alpha_2} \tilde{\cB}
\end{equation}
is a $^*$-isomorphism and that it intertwines $\alpha_1$\ and $\alpha_1'=(\id\ot\id\ot \gamma_1)\mid_{(\cD_1\boxalg^{\alpha_1}\cB)\boxalg^{\alpha_2} \tilde{\cB}}$.
This concludes the proof.

\end{proof}

\begin{remark} Theorem \ref{hierse}\ and proposition \ref{omgekeerd} show that the assignment $\cD_1\recht \cD_1\boxalg^{\alpha_1}\cB$\ and $\alpha_1\recht \id\boxtimes \sde_2$\ yields a bijective  correspondence between actions of $\cG_1$\ and actions of $\cG_2$ (up to conjugacy).
\end{remark}

\subsection{Special cases}

\subsubsection*{The case where $D_1$\ is a homogeneous space}

Suppose we are given two monoidally equivalent compact quantum
groups $\cG_1$ and $\cG_2$ with monoidal equivalence
$\varphi:\cG_2\recht\cG_1$. Denote by $\sde_1$ and $\sde_2$ again
the ergodic actions of full quantum multiplicity on the
corresponding link-algebra $\cB$\ as in formula $(\ref{qqsszz})$.
Consider a quantum subgroup $\cH_1$ of $\cG_1$ and the corresponding
action
\[\de_{\cH_1\backslash\cG_1}:\cC(\cH_1\backslash\cG_1)\to
\cC(\cH_1\backslash\cG_1)\ot \cC(\cG_1)\] on the homogeneous space.

To a homogeneous space of a compact quantum group, naturally there
corresponds a homogeneous space of the link algebra.

\begin{definition}\label{homlinkalg}
We define the homogeneous space $\cB^{\cH_1}$ by
\[\cB^{\cH_1}:=\{a\in\cB\mid (r_{\cH_1}\ot\id)\sde_{1}(a)=1\ot a\}\;.\]
\end{definition}

Note that $\delta_2$\ is an action on $\cB^{\cH_1}$\ because
$\delta_1$\ and $\delta_2$\ commute.

\begin{proposition}\label{moneqhom}
There is a $^*$-isomorphism between $\cC(\cH_1\backslash\cG_1)\boxalg^{\de_{\cH_1\backslash\cG_1}}\cB$\ and
 $\cB^{\cH_1}$. Moreover, this *-isomorphism
intertwines the action $\id\boxtimes \sde_2$\ with the restriction of $\sde_2$
to $\cB^{\cH_1}$.
\end{proposition}
\begin{proof}

 We prove, as in
proposition \ref{link algebra}, that $\sde_1:\cB^{\cH_1}\to \cC(\cH_1\backslash\cG_1)\boxalg^{\de_{\cH_1\backslash\cG_1}}\cB$
is a *-isomorphism. Because
\begin{align*}
(r_{\cH_1}\ot\id)(\de_{\cH_1\backslash\cG_1}\ot\id)\sde_1(a)&=(r_{\cH_1}\ot\id)(\id\ot\sde_1)\sde_1(a)
\\&=(\id\ot\sde_1)(1\ot a)=1\ot\sde_1(a)\;,
\end{align*}
 we get that $\sde_1(\cB^{\cH_1})\subseteq \cC(\cH_1\backslash\cG_1)\boxalg^{\de_{\cH_1\backslash\cG_1}}\cB$.
The injectivity of $\sde_1$  is clear. The surjectivity follows
from the fact that for $a\in \cC(\cH_1\backslash\cG_1)\boxalg^{\de_{\cH_1\backslash\cG_1}}\cB$,
\begin{equation}\label{surjectief}
\sde_1(\eps\ot\id)(a)=(\eps\ot\id\ot\id)(\id\ot\sde_1)(a)=(\eps\ot\id\ot\id)(\de_{\cH_1\backslash\cG_1}\ot\id)(a)=a
\end{equation}
and
\[(r_{\cH_1}\ot\id)\sde_1(\eps\ot\id)(a)=(r_{\cH_1}\ot\id)(a)=(\eps\ot\id)(a),\]
where in the last step we used lemma \ref{coeenheid}.
Because $\sde_1$ and $\sde_2$ commute, $\sde_1$ intertwines
$\id\boxtimes \sde_2$ with the restriction of $\sde_2$ to $\cB^{\cH_1}$. This
ends the proof.
\end{proof}

\subsubsection*{The case where $\alpha_1$\ is the adjoint action $\alpha_{\cG_1}$}
We now look at the special case where we are dealing with the adjoint action, introduced in section $5$:
$$\al_{\cG_1}: \plone\recht\plone\otvn L^{\infty}(\cG_1):\al_{\cG_1}(a) =
\mathbb{V}_1(a\ot 1)\mathbb{V}_1^*\;.$$
We get the following proposition:
\begin{proposition}\label{fundin}
The mapping
\begin{equation}\label{vlees}
\rho:\plont\recht \plone\boxtimes^{\alpha_{\cG_1}} B: \rho(b)=\cX(b\ot 1)\cX^*
\end{equation}
is a $^*$-isomorphism. Moreover, this $^*$-isomorphism intertwines the action $\id\boxtimes \sde_2$\ with the adjoint action $\alpha_{\cG_2}$\ on $\plont$. We also have that
$$\rho:\cL(H_x)\recht \cL(H_{\varphi(x)})\boxalg^{\alpha_{\cG_1}}\cB$$
is a $^*$-isomorphism.
\end{proposition}

\begin{proof}
The following calculation
$$(\id\ot\sde_1)(\mathbb{X}(b\ot 1)\mathbb{X}^*)=(\mathbb{V}_1)_{12}\mathbb{X}_{13}( b\ot1\ot 1)\mathbb{X}^*_{13}(\mathbb{V}_1)_{12}^*=(\alpha_{\cG_1}\ot \id)(\mathbb{X}(b\ot 1)\mathbb{X}^*)$$
shows that $\rho:\plont\recht \plone\boxtimes^{\alpha_{\cG_1}} B$\ is an injective $^*$-homomorphism.
For every $a\in \plone\boxtimes^{\alpha_{\cG_1}} B$, we have that
\begin{align*}
(\id\ot\sde_1)(\mathbb{X}^*a\mathbb{X})&=
\mathbb{X}^*_{13}(\mathbb{V}_1)_{12}^*(\id\ot\sde_1)(a)(\mathbb{V}_1)_{12}\mathbb{X}_{13}
=\mathbb{X}^*_{13}(\mathbb{V}_1)_{12}^*(\alpha_{\cG_1}\ot\id)(a)(\mathbb{V}_1)_{12}\mathbb{X}_{13}=
(\mathbb{X}^*a\mathbb{X})_{13}\;.
\end{align*}
Because $\sde_1$\ is ergodic, it follows that $\mathbb{X}^*a\mathbb{X}=y\ot 1$\ with $y\in \plont$. From this; we get that  for every $a\in \plone\boxtimes^{\alpha_{\cG_1}} B$\ there exists an $y\in \plont$\ such that  $a=\cX(y\ot 1)\cX^*$. This proves the surjectivity of $\rho$. Note that $\rho^{-1}$\ is given by $\rho^{-1}(a)=(\id\ot\om)(\cX^*a\cX)$\ for every $a\in  \plone\boxtimes^{\alpha_{\cG_1}} B$.

We prove now that $(\id\boxtimes \sde_2)\circ \rho=(\rho\ot \id)\circ \alpha_{\cG_2}$. This follows from the following calculation
\begin{align*}
(\id\boxtimes \sde_2)(\rho(b))=(\id\ot\sde_2)(\cX(b\ot 1)\cX^*)=\cX_{12}(V_{2})_{13}(b\ot 1\ot 1)(V_{2})^*_{13}\cX^*_{12}=(\rho\ot \id)(\alpha_{\cG_2}(b))\;.
\end{align*}
It is also immediately clear that $\rho$\ sends $\cL(H_x)$\ to $\cL(H_{\varphi(x)})\boxalg^{\alpha_{\cG_1}}\cB$\ and that $\rho^{-1}$\ sends $\cL(H_{\varphi(x)})\boxalg^{\alpha_{\cG_1}}\cB$\ to $\cL(H_x)$.
This completes the proof.

\end{proof}

\begin{remark} This result suggest strongly that Poisson and Martin boundaries of $\cG_2$\  are related to the bijective construction obtained in this section. In the next two sections, we show that this is indeed the case.
\end{remark}

\section{Poisson boundaries of monoidally equivalent quantum groups}

In this section we prove that the Poisson boundaries of two
monoidally equivalent quantum groups correspond with each other
through the construction of Theorem \ref{hierse}. Recall that, because of remark \ref{vonnie}, we may do all computations immediately on the von Neumann algebraic level.

Consider two monoidally equivalent compact quantum groups $\cG_1$
and $\cG_2$ where the monoidal equivalence is given by $\varphi:
\cG_2\recht \cG_1$\ with corresponding link algebra $B$\ and commuting actions $\sde_1$\ and $\sde_2$.

\begin{notation}\label{not.intertwiners}
From now on, we write respectively
$\Vb_1:=\bigoplus_{x\in\Irred(\cG_1)}U^{\varphi(x)}$\ and
$\Vb_2:=\bigoplus_{x\in\Irred(\cG_2)}U^x$, where $\{
U^{\varphi(x)}\mid\ x\in \Irred(\cG_1)\}$\ and $\{U^{x}\mid \ x\in
\Irred(\cG_2)\}$\ denote the set of irreducible representations of
respectively $\cG_1$\ and $\cG_2$. We also denote by $B=(\cB,\om)''$ the
von Neumann algebraic link algebra of the monoidal equivalence and
by $\Xb:=\bigoplus_{x\in\Irred(\cG_2)}X^x$. We denote the states
$\vphi^1_\mu$ and $\psi^1_\mu$, respectively $\vphi^2_\mu$ and
$\psi^2_\mu$ on $\ell^\infty(\cGh_1)$, respectively
$\ell^\infty(\cGh_2)$.
\end{notation}

Let $\mu$\ be a generating probability measure on $\Irred(\cG_1)$.
Consider the Poisson boundary $\hone$\ of $\cG_1$\ with adjoint
action
\[\alpha_{\cG_1}:\hone\recht \hone\otvn
L^{\infty}(\cG_1):\alpha_{\cG_1}(a)=\mue(a\ot 1)(\mue)^*\;\]
We get the following theorem:

\begin{theorem}\label{eigenspul}
Consider two monoidally equivalent compact quantum groups $\cG_1$
and $\cG_2$\ and let $\mu$\ be a generating probability measure on
$\Irred(\cG_1)$. Then the following
$$\rho: \hont\recht \hone\boxtimes^{\alpha_{\cG_1}}B:\rho(b)=\cX(b\ot 1)\cX^*$$
is a $^*$-isomorphism. Moreover, this $^*$-isomorphism
intertwines the action $\id\boxtimes \sde_2$\ with the adjoint action
$\alpha_{\cG_2}$.
\end{theorem}

\begin{proof}
In the next proposition, we prove that $(P_{\mu,1}\ot\id)\circ\rho=\rho\circ P_{\mu,2}$. Now, because $\rho: \plont\recht \plone\boxtimes^{\alpha_{\cG_1}} B$\ is a $^*$-isomorphism and because of the definition of the product on $\hone$\ and $\hont$,
 we get that $\rho:  \hont\recht \hone\boxtimes^{\alpha_{\cG_1}} B$\ is a $^*$-isomorphism.
 \end{proof}

\begin{proposition} \label{nииии}
Let $\mu$ be a
pro\-bability measure on $\mbox{Irred}(\cG_1)$.  We have that
\begin{equation}\label{zimzam}
(P_{\mu,1}\ot\id)(\rho(b))=\rho( P_{\mu,2}(b))
\end{equation}
for every $b\in \plont$.
\end{proposition}

\begin{proof}
Let $b\in \plont$. We claim that for $x,y\in \Irred(\cG_1)$,
$$(p_x\ot p_y\ot 1)(\deh_1\ot\id)(\rho(b))=X_{13}^xX_{23}^y(\widehat{\Delta}_2(b)\ot
1)\bigr)(X_{23}^y)^*(X_{13}^x)^* \;.$$

Take now $z\in\Irred(\cG_1)$ and $\varphi(T)\in\Mor(x\ot y,z)$. Then
\begin{align}
&(p_x\ot p_y\ot 1)(\deh_1\ot\id)(\rho(b))(\varphi(T)\ot 1) =(\varphi(T)\ot 1)(\rho(b))
=(\varphi(T)\ot 1)\bigl(X^z(b\ot 1)(X^z)^*\bigr)
\notag
\\ &= X_{13}^xX_{23}^y(T\ot 1)(b\ot 1)(X^z)^*=X_{13}^xX_{23}^y(\deh_2\ot\id)(b)(T\ot 1)(X^z)^*\label{iz1}\\
&=X_{13}^xX_{23}^y(\widehat{\Delta}_2(b)\ot
1)\bigr)(X_{23}^y)^*(X_{13}^x)^*(\varphi(T)\ot 1)\notag
\end{align}
where \eqref{iz1} is valid because \[(T\ot
1)(X^z)^*=(X_{23}^y)^*(X_{13}^x)^*(\vphi(T)\ot 1)\;.\] Then, we
get that
\begin{align*}
&(p_x\ot 1)(P_{y,1}\ot \id)(\rho(b))=X^x(\id\ot \psi^1_y\ot
\id)(X_{23}^y(\deh_2(b)\ot 1)(X^y_{23})^*)(X^x)^*\notag\;.
 \end{align*}

We prove that
\begin{equation}\label{beterr}
(\psi_y^1\ot\id)((X^y)(d\ot 1)(X^y)^*)=\psi_y^2(d)1
\end{equation}
for every $d\in \plont$.

If this last equality (\ref{beterr}) is valid, then we get
\begin{align}
(p_x\ot 1)(P_{y,1}\ot\id)(\rho(b))=X^x((p_x\ot 1)(P_{y,2}(b)\ot 1))(X^x)^*=(p_x\ot 1)\rho(P_{y,2}(b))\;.
\end{align}
This means that $(P_{y,1}\ot\id)(\rho(b))=\rho( P_{\mu,2}(b))$\ for every $y\in \Irred(\cG)$\ and (\ref{zimzam}) is true.

The only thing left to prove is formula (\ref{beterr}). Using the definition of $\psi^1_y$ (\ref{vallavalla}), we get
\begin{align}
(\psi_y^1\ot\id)((X^y)(d\ot 1)(X^y)^*)&=(\eta\ot  h_1\ot \id)(U^{\varphi(y)}_{12}X_{13}^y(d\ot 1\ot 1)(X^y_{13})^*(U^{\varphi(y)}_{12})^*)\label{pmutte}
\end{align}
with $\eta$\ any normal state on $\plone$. Then
\begin{align*}
(\ref{pmutte})&= (\eta\ot  h_1\ot \id) ((\id\ot\sde_1)(X^y(d\ot 1)(X^y)^*))=(\eta\ot \om)(X^y(d\ot 1)(X^y)^*)).1
=\psi^2_y(d)1
\end{align*}
where in the last step, we used the orthogonality relations (\ref{eq.orthogonalitymon}). This proves equality (\ref{beterr}) and the proof is complete.
\end{proof}

\section{Applications to Tomatsu's work on Poisson boundaries}
The fundamental result obtained in the last section combined with recent work of Tomatsu on Poisson boundaries makes it possible to identify the Poisson boundary of a large class of quantum groups. This will be the content of this section. First, it provides us with an identification of the Poisson boundary of all duals of compact quantum groups that are monoidally equivalent with a $q$-deformations of a compact Lie group. Moreover, it will enables us to give concrete identifications of Poisson boundaries of some classes of non-amenable discrete quantum groups. Observe that monoidal equivalence does not preserve coamenability. We study in detail the Poisson boundary of the duals of the quantum automorphism groups $A_{aut}(D,\om)$. When $\dim(D)\geq 5$, these are actually non-amenable.

\subsection{Tomatsu's work on Poisson boundaries}

 In \cite{tomatsu1}, Tomatsu has proven that the Poisson
boundary of the dual of a coamenable compact quantum group with
commutative fusion rules can be identified with the homogeneous
space coming from its canonical Kac subgroup (see definition
\ref{Kaccc}).

The main result is the following.

\begin{theorem}[Theorem 4.8 in \cite{tomatsu1}]\label{tom}
Let $\cG$ be a coamenable compact quantum group with commutative
fusion rules and $\cH$ its canonical Kac subgroup. Let $\mu$ be a
generating measure on $\Irred(\cG)$. The Izumi operator
\[\Phi: L^{\infty}(\cH\backslash\cG)\to
H^{\infty}(\cGh,\mu):a\mapsto (\id\ot h)(\Vb^*(1\ot a)\Vb)\] is a
*-isomorphism and intertwines the adjoint action $\alpha_{\cG}$\
with the action $\Delta_{\cH\backslash\cG}$\ as defined in remark
$\ref{qawz}$.
\end{theorem}

\begin{remark} The mapping $\Phi$\ is the Izumi operator, introduced in \cite{iz1}.
\end{remark}
 This
gives us immediately the Poisson boundary of a whole class of
quantum groups. Moreover, for $q$-deformations of classical
Lie-groups, Tomatsu proves that the canonical Kac subgroup is just the maximal torus.

\subsection{Identification of a large class of Poisson boundaries}

Using the previous section, we obtain the Poisson boundary of
every compact quantum group with commutative fusion rules which is
monoidally equivalent to a coamenable one. Moreover, we obtain a
concrete description of the Poisson boundary as a homogeneous
space of the link algebra.

So, consider a compact quantum group $(\cG_1,\de_1)$ which is
coamenable and has commutative fusion rules. Let $(\cG_2,\de_2)$
be monoidally equivalent with $(\cG_1,\de_1)$\
with monoidal equivalence given by $\vphi:\cG_2\to\cG_1$. Again,
we have the link algebra $B$\ and the two commuting actions $\delta_1$\ and $\delta_2$\ as before. Denote by
$\cH_1$ the canonical Kac group of $\cG_1$ and by
$r_{\cH_1}:\cC(\cG_1)\to \cC(\cH_1)$\ be the corresponding
restriction map.

Tomatsu's result combined with theorem \ref{eigenspul} and proposition \ref{moneqhom}
  gives us the following theorem:

\begin{theorem}\label{mon.amenable}
Consider a coamenable compact quantum group $\cG_1$ with
commutative fusion rules. Let $\cG_2$ be a compact quantum groups
that is monoidally equivalent to $\cG_1$. Denote by $B$ the von
Neumann algebraic link algebra associated to the monoidal
equivalence. Let $\cH_1$ be the canonical Kac subgroup of $\cG_1$.
Consider a generating measure $\mu$ on $\Irred(\cG_1)$. Then
$(\hont,\alpha_{\cG_2})$ is isomorphic to $(B^{\cH_1},\sde_2)$.
The isomorphism is given by the following generalized Izumi
operator
\begin{equation}\label{izzzzz}
\Theta:B^{\cH_1}\recht \hont:a\mapsto
(\id\ot\om)(\Xb^*(1\ot a)\Xb)\;.
\end{equation}
This $^*$-isomorphism intertwines the adjoint action $\alpha_{\cG_2}$
and the action $\sde_2$.

\end{theorem}

\begin{proof}
By Tomatsu, $\Phi:L^\infty(\cH_1\backslash\cG_1)\to \hone$ is a
$^*$-isomorphism. Because $\Phi$ intertwines the actions
$\de_{\cH_1\backslash\cG_1}$ and $\alpha_{\cG_1}$,
\[\Phi\ot\id:L^{\infty}(\cH_1\backslash\cG_1)\boxtimes^{\de_{\cH_1\backslash\cG_1}} B\to \hone\boxtimes^{\alpha_{\cG_1}} B\] is also a $^*$-isomorphism. Combining this with Theorem \ref{eigenspul}, it
follows that $L^{\infty}(\cH_1\backslash\cG_1)\boxtimes^{\de_{\cH_1\backslash\cG_1}} B$ and $\hont$ are $^*$-isomorphic through
$\rho^{-1}\circ(\Phi\ot\id)$.
It follows from proposition \ref{moneqhom} that $\sde_1:B^{\cH_1}\to
L^{\infty}(\cH_1\backslash\cG_1)\boxtimes^{\de_{\cH_1\backslash\cG_1}} B$ is a $^*$-isomorphism. Hence
\[\rho^{-1}\circ(\Phi\ot\id)\circ\sde_1:B^{\cH_1}\to \hont\] is a *-isomorphism.

Now we just need to prove that
$\rho^{-1}\circ(\Phi\ot\id)\circ\sde_1=\Theta$, which follows
from the next obvious calculation.
\begin{align*}
\rho^{-1}(\Phi\ot\id)\sde_1(a)&=(\id\ot\om)(\Xb^*(\id\ot
h_1\ot\id)((\mue)^*_{12}(\sde_1(a))_{23}(\mue)_{12})\Xb)
\\&=(\id\ot\om)(\id\ot h_1\ot\id)(\id\ot\sde_1)(\Xb^*(1\ot a)\Xb)\\
&=(\id\ot\om)\bigl(((\id\ot\om)(\Xb^*(1\ot a)\Xb))\ot 1\bigr)
=\Theta(a)
\end{align*}
for all $a\in B^{\cH_1}$.

Moreover, $\Theta$ intertwines $\alpha_{\cG_2}$ and
$\sde_2$\ because for all $a\in B^{\cH_1}$,

\begin{align*}
(\Theta\ot\id)\sde_2(a)&=(\id\ot\om\ot\id)(\Xb^*_{12}\sde_2(a)_{23}\Xb_{12})
\\&=(\id\ot\om\ot\id)((\mut)_{13}(\id\ot\sde_2)(\Xb^*(1\ot a)\Xb)(\mut)^*_{13})\\
&=\mut((\id\ot\om)(\Xb^*(1\ot a)\Xb)\ot 1)(\mut)^*
=\alpha_{\cG_2}(\Theta(a))\;.
\end{align*}

This completes the proof.

\end{proof}

Now we have enough material to identify the Poisson boundary of some classes of discrete quantum groups which are not amenable. A first important class of quantum groups that satisfy this are the duals of $A_o(F)$. If $\dim(F)\geq 3$, then $A_o(F)$\ is not coamenable.
The Poisson boundary of their dual was already obtained in a different way
(but also using monoidal equivalence) by Vaes and the second author in \cite{VaeVa}. In fact, they started by constructing a generalized Izumi operator
as in formula $(\ref{izzzzz})$\ for the specific case of $A_o(F)$.
They proved that this Izumi operator is multiplicative on
$L^{\infty}(A_o(F,F_q))^{\mathbb{T}}$\ by using the monoidal
equivalence of $A_o(F)$ and $\SU_q(2)$. Hence, they reduced the
identification problem to a purely $\SU_q(2)$-problem. However, as
every $A_o(F)$ is monoidally equivalent to some $SU_q(2)$, we can
identify the Poisson boundary of $\widehat{A_o(F)}$ also using
theorem \ref{mon.amenable}, which is much more general.

Another, new class of examples will come from quantum automorphism groups $A_{aut}(D,\om)$, which we
explore in the next section.

\subsection{Examples: Quantum automorphism groups.}

In this section we identify the Poisson
boundary for  $\widehat{A_{aut}(D,\om)}$ with $D$ a \cst-algebra
of finite dimension $\geq4$ and $\om$ a $\sde$-state on $D$. To
do this, we make use of the previous section.

From theorem \ref{extrabew}, it follows that every quantum
automorphism group of this type is monoidally equivalent with one
of the form $A_{aut}(M_2(\C),\Tr(\cdot F))$, where
$\Tr(F^{-1})=\sde^2$. Because of the quantum Kesten result (see
\cite{banica7}), it follows that $A_{aut}(M_2(\C),\Tr(\cdot F))$
is coamenable. Moreover, it has the fusion rules of $SO(3)$ and
those are commutative. Hence, we can apply theorem \ref{tom} of
Tomatsu. We now prove that the canonical Kac subgroup is just the one-dimensional
torus $\mathbb{T}$.

\subsubsection*{The canonical Kac subgroup of
$A_{aut}(M_2(\C),\Tr(\cdot D))$}

Denote by $\cG$ the compact quantum group $A_{aut}(M_2(\C),\Tr(\cdot
D))$ and by $\cH$ its canonical Kac subgroup. Observe that we can
take $D$ a diagonal matrix. We consider only non-trivial $D$ here. If $D$ equals the identity matrix,
we just get the compact group $\SO(3)$ which is already Kac. Hence $\widehat{\SO(3)}$ has trivial Poisson boundary. Denote by $\pi:C(\cG)\to C(\cH)$ the
canonical projection map. We denote by $U$ the fundamental
irreducible representation with label $1$ and by $Q$ the matrix
corresponding to $U$ as defined in $\ref{eq.adjointmap}$, normalized
such that $\Tr(Q)=\Tr(Q^{-1})$. The eigenvalues of $Q$ are of the
form $1,q,q^{-1}$.

Now $V:=(\id\ot\pi)(U)$ is a representation of $\cH$ and because
$\cH$ is Kac, $\overline{V}=(\id\ot\pi)(\overline{U})$ must be
unitary. The matrix $F=\sqrt{Q^T}$, unitarizes $\overline{U}$, what
in this case means $(F\ot 1)\overline{U}(F^{-1}\ot 1)=U$. We claim
that $V$ breaks up in $3$ one-dimensional representations. As every
representation of $\cH$ appears in a repeated tensor power of $V$,
it follows that all irreducible representations of $\cH$ have
dimension one.

\emph{Proof of claim:} As $F^*F$ has $3$ different eigenvalues, it
suffices to prove that $F^*F$ and $V$ commute. It holds that
$V=F\overline{V}F^{-1}$, so
\[VF=F\overline{V}\quad\text{and}\quad F^*V^*=\overline{V}^*F^*\;.\]
As $\overline{V}$ is unitary, it follows that $VFF^*V^*=FF^*$, which
means that $F^*F\in\End(V)$.

Since all irreducible representations of $\cH$ have dimension $1$,
we conclude that $\cH$ is the dual of a discrete group $\Gamma$.
Denote by $u_g$ the irreducible representation of $\cH$
corresponding to $g\in\Gamma$. Since $V\cong\overline{V}$, there
exist $g,h\in\Gamma$ such that
\[V \cong u_g \oplus u_h \oplus u_{g^{-1}}\;.\]
Observe that $\Gamma$ is generated by $g$ and $h$. We claim that
$\Gamma$ is abelian. Since $U$ is a subrepresentation of $U^{\ot
2}$, $V$ is a subrepresentation of $V^{\ot 2}$. But
\[V^{\ot 2}=u_{g^2}\oplus u_{h^2}\oplus u_{g^{-2}}\oplus
2u_e\oplus u_{hg}\oplus u_{gh}\oplus u_{hg^{-1}}\oplus
u_{g^{-1}h}\;,\] implying that $h\in\{g^2,h^2, g^{-2},e, hg,
gh,hg^{-1},g^{-1}h\}$. Any of these possibilities for $h$ imply that
the group generated by $g$ and $h$ is abelian.

As $\Gamma$ is commutative, $\cH$ is just a commutative compact
group. Hence the maximal quantum subgroup of Kac type of $\cG$ is
the maximal compact subgroup of $\cG$.

Suppose $\chi:C(\cG)\to\C$ is a character and $\alpha:M_2(\C)\to
M_2(\C)\ot C(\cG)$ the canonical action of $\cG$ on $M_2(\C)$ coming
from $U$. Now $(\id\ot\chi)\alpha$ is an automorphism of $M_2(\C)$,
and hence implemented by a unitary matrix $A$. Moreover, as
$\Tr(\cdot D)$ is invariant under $(\id\ot\chi)\alpha$,
\[\Tr(DAxA^*)=\Tr(Dx)\quad\text{ for all}\quad x\in M_2(\C)\;.\]
Hence $A$ is a diagonal matrix. But then
$\Ad(A)=\Ad(\diag(z,\bar{z}))$ for some $z\in \mathbb{T}$.

On the other hand, $\mathbb{T}$ acts on $M_2(\C)$ by $\Ad(\diag(z,\bar{z}))$. This action $\sde$ is $\Tr(\cdot D)$-invariant, so
because of the universality of $\cG$, there exists a morphism of
quantum groups $\pi: C(\cG)\to C(\mathbb{T})$ such that
$(\id\ot\pi)\alpha=\sde$.

We may conclude that the maximal Kac subgroup of $\cG$ is the one-dimensional torus
$\mathbb{T}$.

\subsubsection*{The Poisson boundary of $A_{aut}(D,\om)$}

Theorem \ref{extrabew} gives a nice description of the link algebra of two monoidally equivalent quantum automorphism groups $A_{aut}(D_1,\om_1)$\ and $A_{aut}(D_2,\om_2)$. Together with theorem \ref{mon.amenable} we obtain the following
result:

\begin{theorem}
Consider $A_{aut}(D,\om)$ with $D$ a \cst-algebra of finite
dimension strictly bigger than $4$ and $\vphi$ a $\sde$-state on $B$. Take
$F\in M_2(\C)$ such that $\Tr(F^{-1})=\sde^2$. Then the Poisson
boundary of $\widehat{A_{aut}(D,\om)}$ is given by
$L^\infty(A_{aut}((D,\om),(M_2(\C),\Tr(\cdot F))))^\mathbb{T}$.
\end{theorem}

The $4$-dimensional case was considered above, except for the case $A_{aut}(\C^4)$. But this compact quantum group is coamenable by the quantum Kesten result and moreover Kac, so $\widehat{A_{aut}(\C^4)}$ has trivial Poisson boundary. This completes the identification of Poisson boundaries of the duals of quantum automorphism groups.

\section{The Martin boundary of monoidally equivalent quantum groups}
We prove that the Martin boundaries of the duals of two monoidally equivalent compact quantum groups are related to each other through the construction of theorem \ref{hierse}.

So, again, we start from a monoidal equivalence $\varphi:\cG_2\recht\cG_1$\ with link algebra $B$\ and commuting actions $\sde_1$\ and $\sde_2$.

\begin{theorem}\label{eigenstuk2}
Consider two monoidally equivalent compact quantum groups $\cG_1$
and $\cG_2$ and let $\mu$\ be a generating probability measure on
$\Irred(\cG_1)$. Suppose $A_{\mu,1}$\ is the Martin boundary of the discrete quantum group $\cGh_1$. Then $$\rho:\widetilde{A}_{\mu,2}\recht\widetilde{A}_{\mu,1}\boxtimes_{red}^{\alpha_{\cG_1}}B:b\mapsto \cX(b\ot 1) \cX^*$$ are also
 $^*$-isomorphisms. This $^*$-isomorphism intertwines the
action $\id\boxtimes \sde_2$\ with the action $\alpha_{\cG_2}$. Moreover, we have that $\rho:c_0(\cGh_2)\recht c_0(\cGh_1)\boxtimes_{red}^{\alpha_{\cG_1}}B $\ is a $^*$-isomorphism, by which
$$\rho:A_{\mu,2}\recht A_{\mu,1}\boxtimes_{red}^{\pi_{\cG_1}}B:b\mapsto \cX(b\ot 1) \cX^*$$
is also a $^*$-isomorphism which intertwines the
action $\id\boxtimes \sde_2$\ with the action $\pi_{\cG_2}$.
\end{theorem}

\begin{proof}
 Recall that
$$\rho:\cL(H_x)\recht \cL(H_{\varphi(x)})\boxalg^{\alpha_{\cG_1}} \cB$$
is a $^*$-isomorphism. So it follows that
$$\rho: c_c(\cGh_2) \recht c_c(\cGh_1)\boxalg^{\alpha_{\cG_1}} \cB\quad \mbox{and}\quad \rho:c_0(\cGh_2)\recht c_0(\cGh_1)\boxtimes_{red}^{\alpha_{\cG_1}} B $$
is also a $^*$-isomorphism.
Proposition \ref{nииии} gives us that
$$(P_{\omu,1}\ot\id)(\rho(b))=\rho( P_{\omu,2}(b))\quad \mbox{and thus}\quad (G_{\omu,1}\ot\id)(\rho(b))=\rho( G_{\omu,2}(b))$$\ for every $b\in c_c(\cG_2)$.  Because $G_{\omu}(p_\varepsilon)\in
M(c_c(\Irred(\cG)))$\ for every $\mu\in \Irred(\cG)$, we get

\begin{align*}
\rho( K_{\omu,2}(b))=\cX(G_{\omu,2}(p_\varepsilon)^{-1}G_{\omu,2}(b)\ot 1)\cX^*
=(G_{\omu,1}(p_\varepsilon)^{-1}\ot 1)\rho(G_{\omu,2}(b))
=(K_{\omu,1}\ot\id)(\rho(b))\;.
\end{align*}
This fact, combined with the fact that $\rho:\plont\recht \plone\boxtimes^{\alpha_{\cG_1}} B$\ is a $^*$-isomorphism, gives that
 $$\rho: \widetilde{A}_{\mu,2}\recht \widetilde{A}_{\mu,1}\boxtimes_{red}^{\alpha_{\cG_1}}B$$
is an injective $^*$-homomorphism.

The only thing left to prove is that $\rho$\ is surjective. We consider the inverse monoidal equivalence
$\varphi^{-1}:\cG_1\recht \cG_2$.  So, in the same way as above, we get that
 $$\widetilde{\rho}: \widetilde{A}_{\mu,1}\recht \widetilde{A}_{\mu,2}\boxtimes_{red}^{\alpha_{\cG_2}}\widetilde{B}:\widetilde{\rho}(b)=\mathbb{Y}(b\ot 1)\mathbb{Y}^*\; ,$$
is an injective $^*$-homomorphism. But we have that
$$\Lambda:(\widetilde{\rho}\ot \id)\circ\rho:\widetilde{A}_{\mu,2}\recht (\widetilde{A}_{\mu,2}\boxtimes_{red}^{\alpha_{\cG_2}}\tilde{B})\boxtimes_{red}^{\alpha_{\cG_1}} B:\Lambda(b)=
 \Yb_{12}\Xb_{13}(b\ot 1\ot 1)\Xb_{13}^*\Yb^*_{12}$$
is just the natural $^*$-isomorphism from proposition \ref{omgekeerd} for the monoidal equivalence $\varphi^{-1}:\cG_1\recht \cG_2$. Indeed, $\Lambda=(\id\ot\pi)\alpha_{\cG_2}$ with \[\pi:L^\infty(\cG_2)\to \widetilde{B}\boxtimes^{\gamma_1}B\quad\text{ given by }\quad (\id\ot\pi)(\Vb_2)=\Yb_{12}\Xb_{13}\;.\] So,
$\rho:\widetilde{A}_{\mu,2}\recht\widetilde{A}_{\mu,1}\boxtimes_{red}^{\alpha_{\cG_1}}B$\ is also surjective and thus a $^*$-isomorphism.
  In proposition \ref{fundin}, we showed that $\rho$\  intertwines the
action $\id\boxtimes \sde_2$\ with the action $\alpha_{\cG_2}$. We have also shown that $\rho:c_0(\cGh_2)\recht c_0(\cGh_1)\boxtimes_{red}^{\alpha_{\cG_1}}B $\ is a $^*$-isomorphism, such that
$$\rho: A_{\mu,2}\recht\mathcal{A}_{\mu,1}\boxtimes_{red}^{\alpha_{\cG_1}}B \; ,$$
is a $^*$-isomorphism which  intertwines the
actions $\id\boxtimes \sde_2$\ and $\pi_{\cG_2}$. This completes the proof of the theorem.

\subsection*{Examples: The universal orthogonal quantum groups $A_o(F)$}
Given an $A_o(F)$, there exists a $q\in]-1,1[\ \{0\}$\ such that $A_o(F)\meq SU_q(2)$.
In \cite{NT}, Neshveyev and Tuset identified the Martin boundary of $\widetilde{SU_q(2)}$\ ($0<\vert q\vert<1$), under the restriction that the measure $\mu$\ has finite first moment, with the Podl\`{e}s sphere $C(\mathbb{T}\backslash SU_q(2))$. By theorem $\ref{eigenstuk2}$, we get that $A_{\mu,2}\cong C(\mathbb{T}\backslash SU_q(2))\boxtimes^{\pi_{\cG_1}}_{red}\cB$. But $C(\mathbb{T}\backslash SU_q(2))$\ is a quotient space of $L^{\infty}(\cG_1)$, so in exactly the same way as in theorem $\ref{mon.amenable}$, we show that $$C(\mathbb{T}\backslash SU_q(2))\boxtimes^{\pi_{\cG_1}}_{red}\cB\cong B_{r}^{\mathbb{T}}$$
where $B_{r}^{\mathbb{T}}$\ is the norm-closure of $\cB^{\mathbb{T}}\subseteq B^{\mathbb{T}}$\ and $\cB$\ the link algebra. Theorem
$\ref{constructionk}$\ says that $\cB\cong \mathcal{C}(A_o(F,F_q))$.
Combining these facts, we obtain the following theorem:

\begin{theorem} \label{thm.martinbeter}
Let $\cG=A_o(F)$. Then there exists  $q\in ]-1,1[\setminus \{0\}$\ such that $\cG\meq SU_q(2)$. Take  $\mu$\ a generating measure on $\Irred(\cG)$ that is transient
and has finite first moment:
$$\sum_{x \in \N} x \mu(x) < \infty \; .$$
Then the Martin boundary $A_{\mu}$\ of $\cG$\ is $^*$-isomorphic with  $(C(A_o(F,F_q)))_{r}^{\mathbb{T}}$.
  Moreover, this $^*$-isomorphism
intertwines the action $\sde_2$\ with the action
$\pi_{\cG_2}$.
\end{theorem}

\end{proof}

\begin{remark}
The above identification was already obtained by Vaes and the second author in \cite{VaeVa} by using another method supporting on techniques from \cite{Vbo} which allow, in the case of $A_o(F)$, to deduce the Martin boundary from the Poisson boundary. The result by which we obtain the identification here is more direct and much more general.
\end{remark}

\end{document}